\documentclass[12pt]{article}

\usepackage{amssymb}
\usepackage{amsmath}
\usepackage{comment}
\usepackage{latexsym}
\usepackage{enumerate}
\usepackage{hyperref}
\newtheorem{theorem}{Theorem}[section]
\newtheorem{proposition}[theorem]{Proposition}
\newtheorem{corollary}[theorem]{Corollary}

\newtheorem{remark}[theorem]{Remark}

\newtheorem{definition}[theorem]{Definition} 
\newcommand{\half}{\frac{1}{2}}
\newcommand{\tpi}{2\pi i}
\newcommand{\der}[1]{\ensuremath{\partial_{#1}}}
\newcommand{\Z}{\mathbb{Z}}
\newcommand{\C}{\mathbb{C}}
\newcommand{\R}{\mathbb{R}}
\newcommand{\HH}{\mathbb{H}}
\newcommand{\vac}{\mathbf {1}}
\newcommand{\delz}{\partial_{z}} 
\newcommand{\wtil}{\ensuremath{\tilde{\omega}}}
\newcommand{\Drho}{{\mathcal D}^{\rho}} 
\newcommand{\Mcal}{{\mathcal M}} 
\newcommand{\Ycal}{{\mathcal Y}}
\newcommand{\qhat}{{\widehat q}}
\newcommand{\abar}{\overline{a}}

\newcommand{\xhat}{\widehat{x}}

\DeclareMathOperator{\diag}{diag}
\DeclareMathOperator{\SL}{SL}
\DeclareMathOperator{\Sp}{Sp}
\DeclareMathOperator{\Tr}{Tr}
\DeclareMathOperator{\wt}{wt}
\DeclareMathOperator{\Id}{Id}

\DeclareMathOperator{\Lin}{Lin}
\begin{document}
\normalfont
\title{Genus Two Partition and Correlation Functions for Fermionic Vertex Operator Superalgebras~II} 
\author{
Michael P. Tuite
 and Alexander Zuevsky\thanks{
Supported by a Science Foundation Ireland Research Frontiers Programme Grant, and
by Max--Planck Institut f\"{u}r Mathematik, Bonn}
\\
School of Mathematics, Statistics and Applied Mathematics, \\
National University of Ireland Galway \\
University Road, Galway, Ireland}
\maketitle
\begin{abstract}
We define and compute the continuous orbifold partition function and a generating function for all $n$-point correlation functions for the rank two free fermion vertex operator superalgebra on a genus two Riemann surface formed by self-sewing a torus.
The partition function is proportional to an infinite dimensional determinant
with entries arising from torus Szeg\H{o} kernel and the generating function is proportional to a finite determinant of genus two Szeg\H{o} kernels. 
These results follow from an explicit analysis of all torus $n$-point correlation functions for intertwiners of the irreducible modules of the Heisenberg vertex operator algebra.
We prove that the partition and $n$-point correlation functions
are holomorphic on a suitable domain and describe their modular properties. We also describe an identity for the genus two Riemann theta series analogous to the Jacobi triple product identity.
\end{abstract}
\maketitle
\newpage

\section{Introduction}
This paper is a continuation of the study   in the 
 theory of Vertex Operator Algebras (VOAs)  \cite{FLM, FHL, K, MN, MT5}
of partition functions and correlation functions on a genus two Riemann surface \cite{T, MT1, MT2, MT3, MT4, TZ2}.
Such a surface can be constructed by sewing two separate tori together, which we refer to as the 
$\epsilon$--formalism, or by self-sewing a single torus, which we refer to as 
the $\rho$--formalism \cite{MT2}. 
In our approach, we define the genus two partition and $n$--point correlation functions in terms of genus one VOA data combined in accordance with the given sewing scheme. 	
This is based solely on the properties of the  VOA with
no assumed analytic or modular properties for partition or correlation functions.
This is in contrast to other approaches to correlation functions using algebraic geometric methods in mathematics (e.g. \cite{TUY, Z2}) or other approaches taken in physics (e.g. \cite{EO,  DP, DVPFHLS}).

In the companion paper ref.~\cite{TZ2} we studied the continuous orbifold partition and correlation functions for the rank two free fermion Vertex Operator Superalgebra (VOSA) $V(H,\Z+\half )^{\otimes 2}$ on a genus two surface  in the $\epsilon$-formalism. 
The present paper studies this system in the $\rho$--formalism where a genus two Riemann surface is constructed by sewing a handle  at two points, separated by $w$, to a torus with modular parameter $\tau$ and with sewing parameter $\rho$ giving a family of surfaces parametrized by $(\tau ,w,\rho )\in  \Drho$, a particular sewing domain. 
Although there are some features in common with ref.~\cite{TZ2}, the present $\rho$--formalism approach is more subtle in two important respects:
firstly, we must consider genus one correlation functions for intertwining vertex operators associated with $g$-twisted modules for a continuous automorphism $g$ and secondly,  the genus two partition and $n$--point functions obtained are multi-valued on $\Drho$.
There is no Zhu reduction theory available for general genus one intertwiner correlation functions unlike the the reduction theory developed in \cite{MTZ} and used in the $\epsilon$-formalism.  
Instead, we use bosonization and compute all torus  $n$--point functions for the Generalized VOA $\Mcal=\oplus_{\alpha\in\C}M\otimes e^{\alpha}$ formed by the intertwiners for the rank one Heisenberg  VOA $M$ and all of its irreducible modules $M\otimes e^{\alpha}$ for all $\alpha\in \C$ \cite{TZ3}. 
This requires an extension of previous results for genus one lattice theory correlation functions to the Generalized VOA $\Mcal$.
We also make repeated use of an explicit expression of genus two Szeg\H{o} kernel in terms of genus one data in the $\rho$--formalism described in \cite{TZ1}.  
The genus two partition and $n$-point  functions are computed explicitly and involve an infinite determinant whose components arise from genus one Szeg\H{o} kernel data. 
The genus two partition and $n$-point  functions  are multi-valued functions on $\Drho$ that lift to holomorphic functions on a suitable covering space and are modular invariant under a subgroup of the genus two symplectic group $\Sp(4,\Z)$.
We also compute the partition function in a bosonic basis and obtain an analogue of the classic Jacobi  triple product formula for the genus two Riemann theta series.

\medskip

Section~\ref{Section_Rho_g} contains a review of the $\rho$-formalism for the construction of a genus two Szeg\H{o} kernel $S^{(2)}$ on $\Drho$ in terms of elliptic Szeg\H{o} kernel data  and given  multipliers on the four cycles of a genus two Riemann surface formed by self-sewing a torus \cite{MT2, TZ1}.
We introduce  a certain infinite matrix $T$ constructed from the genus one Szeg\H{o} data and the  multipliers which plays a central role: $(I-T)^{-1}$ appears explicitly in the formula for $S^{(2)}$ and $\det(I-T)$, which is holomorphic on $\Drho$, appears in our later VOA results.

Section~\ref{vosas} reviews Vertex Operator Superalgebras (VOSA) and $g$-twisted modules with particular emphasis on the rank two free fermion  VOSA $V(H,\Z+\half )^{\otimes 2}$ for a continuous automorphism $g=e^{-2\pi i\alpha a(0)}$ for Heisenberg vector $a$. 
We often make use of bosonization whereby 
$V(H,\Z+\half )^{\otimes 2}\cong V_{\Z}$, the $\Z$--lattice VOSA with trivial cocycle structure.
In particular, we consider the the interwining operators  $\Mcal=\oplus_{\alpha\in\C}M\otimes e^{\alpha}$ that satisfy a Generalized VOA described in detail in \cite{TZ3}. This leads to 
generalized notions of locality, skew-symmetry, associativity, commutivity with respect to the Heisenberg sub VOA and an invariant linear form on $\Mcal$.
 
In Section~\ref{Torus} we compute general $n$-point intertwiner correlation functions for elements of $\Mcal$ in terms of elliptic prime forms and elliptic and quasi-modular forms. These expressions are natural generalizations of results for lattice VOAs \cite{MT1} extended to $\Mcal$. This allows us to compute a generating form for all intertwiner torus $n$--point functions for $V_{\Z}$ expressed in terms of a finite determinant of genus one Szeg\H{o} kernels.

Section~\ref{genus_two_n-point} contains our main results concerning the genus two partition function and a generating form for all genus two $n$--point functions. 
These are defined in terms of four commuting automorphisms generated by $a(0)$ which determine the $S^{(2)}$ multipliers.  
The partition function which is defined as a sum over a basis for a $g$--twisted module of particular genus one intertwining 2--point functions is explicitly computed in Theorem \ref{Theorem_Z2_ferm}. The result is proportional to $\det(I-T)$ but is multivalued on $\Drho$ due to a separate factor of $\left(-\rho/K(w,\tau)^2\right)^{\half\kappa^2}$ where $g=e^{-2\pi i\kappa a(0)}$ for $\kappa\in(-\half,\half)$ and $K(w,\tau)$ is  the elliptic prime form. 
In Theorem~\ref{Theorem_modul_inv}
we show that the partition function lifts to a holomorphic function on a covering space  $\widehat{\mathcal D}^{\rho}$
and prove that it is modular invariant under a particular subgroup $L\subset \Sp(4,\Z)$, the genus two symplectic group.
We also compute the partition function in an alternative bosonic basis and obtain an analogue of the Jacobi  triple product formula for the genus two Riemann theta series in Theorem~\ref{theorem:Triple}.
Lastly, in Theorem~\ref{generating_n_point_rank_two} a modular invariant generating form for all genus two $n$-point correlation functions is  computed in terms of the genus two partition function and a finite determinant of genus two Szeg\H{o} kernels.
\section{The Szeg\H{o} Kernel in the $\rho$-Formalism}
\label{Section_Rho_g} 
In this section we consider the construction
of a genus two Riemann surface formed by self-sewing a handle to a torus using an explicit sewing scheme which we refer to as the $\rho$-formalism.
In particular, we review the construction of an explicit formula for the genus two Szeg\H{o} kernel in terms of Szeg\H{o} kernel on the torus. Further details can be found in   \cite{MT2} and \cite{TZ1}. 
\subsection{The Szeg\H{o} Kernel on a Riemann Surface}  
\label{Szegokernel}
Consider a compact connected Riemann surface $\Sigma$ of genus $g$ with canonical
homology cycle basis $a_{i},b_{i}$ for $i=1,\ldots, g$.
 Let  $\nu_{i}$ be a basis of holomorphic 1-forms 
with normalization $\oint_{a_{i}}\nu_{j}=2\pi i\delta_{ij}$ and
period matrix $\Omega_{ij}=\frac{1}{2\pi i}\oint_{b_{i}}\nu_{j}\in \mathbb{H}_g$, 
the Siegel upper half plane.
Define the theta function with real characteristics  \cite{M, F, FK}
\begin{equation}
\label{theta}
\vartheta 
\begin{bmatrix}
\alpha \\ \beta  
\end{bmatrix}
 \left( {z} \vert   \Omega \right)  
= 
\sum\limits_{ n \in {\Z^{g}} } e^{
 i \pi (n + {\alpha}).\Omega .(n+{\alpha}) + 
 (n+ {\alpha}). ({{z}+ 2 \pi i {\beta}})},
\end{equation}
for $\alpha=(\alpha_j),{\beta}=(\beta_j)\in \mathbb{R}^g$ and ${z}=(z_j)\in \mathbb{C}^g$ for $j=1,\ldots, g$. 
The Szeg\H{o} kernel 
is defined for $\vartheta \left[
\begin{smallmatrix}
\alpha \\ \beta  
\end{smallmatrix} \right]
( 0\vert \Omega)\neq 0$ by \cite{S}, \cite{HS}, \cite{F}
\begin{equation}
\label{gSzego10}
 S
\begin{bmatrix}
\theta  \\ \phi  
\end{bmatrix}
(x, y\vert \Omega) =
\frac
{ \vartheta
\begin{bmatrix}
\alpha \\ \beta  
\end{bmatrix}
\left( \int_{y}^{x}\nu  \,\vert  \Omega\right)} 
 {\vartheta
 \begin{bmatrix}
\alpha \\ \beta  
\end{bmatrix}
( 0\vert \Omega)  E(x, y)},
\end{equation}
where 
$\theta=({\theta}_j),\ \phi=(\phi_j)\in U(1)^n$ for
\begin{equation}
{\theta}_j=-e^{-2 \pi i \beta_j}, \quad{\phi}_j= -e^{2 \pi i \alpha_j},\quad j=1,\ldots, g. 
\label{eq:periodicities}
\end{equation}
(where the $-1$ factors are included for later convenience)
and $E(x,y)$ is the prime form \cite{M, F}. 
The Szeg\H{o} kernel is periodic in
 $x$ along the $a_i$ and $b_j$ cycles 
 with multipliers $-\phi_i$ and $-\theta_j$ respectively and is 
a meromorphic $(\half ,\half )$-form (and is thus necessarily defined on a double-cover $\widetilde{\Sigma}$ of the Riemann surface) satisfying\footnote{We use the convention $E(x,y)\sim (x-y )dx^{-\half }dy^{-\half }$ for $x\sim y$.}
\begin{eqnarray*}
 S\begin{bmatrix}
\theta  \\ \phi  
\end{bmatrix}(x,y) 
& \sim &  \frac{1}{x-y}\; dx^{\half }\; dy^{\half }
\quad \mbox{for } x\sim y,
 \label{Sz_local}\\
 S\begin{bmatrix}
\theta  \\ \phi  
\end{bmatrix}(x,y) 
& = &   -S\begin{bmatrix}
\theta^{-1}  \\ \phi^{-1}  
\end{bmatrix} (y,x),\label{Sz_skewsym}
\end{eqnarray*}
where $\theta^{-1}=({\theta}_i^{-1})$ and $\phi^{-1}=(\phi_i^{-1})$. 
The Szeg\H{o} kernel transforms under the action of the symplectic  modular group
for  $\left[
\begin{smallmatrix}
A & B \\ 
C & D
\end{smallmatrix}
\right]
\in \Sp(2g,\Z)$  as follows \cite{F}
\begin{eqnarray}
S\begin{bmatrix}
\tilde{\theta}  \\  \tilde{\phi}
\end{bmatrix}
(x, y\vert   \widetilde\Omega)
&=&
S\begin{bmatrix}
\theta  \\ \phi  
\end{bmatrix} (x, y\vert \Omega),
\label{Szmod}
\\
\widetilde \Omega
&=&
\left(A\Omega+B\right)\left(C\Omega+D\right)^{-1},
\label{eq:modOmega}
\end{eqnarray} 
 where $\tilde{\theta}_j=-e^{-\tpi  \tilde{\beta}_j}$, $\tilde{\phi}_j= -e^{\tpi  \tilde{\alpha}_j}$ with
\begin{equation}
\begin{bmatrix}
-\tilde{\beta} \\ 
\tilde{\alpha}
\end{bmatrix}
=
\begin{bmatrix}
A & B \\ 
C & D
\end{bmatrix}
\begin{bmatrix}
{-\beta}\\ 
\alpha
\end{bmatrix}
+\half 
\begin{bmatrix}
{-\diag(AB^T)} \\ 
{\diag(CD^T)}
\end{bmatrix},
\label{albetatilde}
\end{equation}
where $\diag(M)$ denotes the diagonal elements of a matrix $M$.
\subsection{Genus Two Riemann Surface Formed from Self-Sewing a Torus}
\label{Sect_Rho_g}
We now briefly review the construction of a genus two Riemann surface 
$\Sigma^{(2)}$ formed by self-sewing a handle to a torus
which we refer to this as the $\rho$-formalism \cite{MT2,TZ1}. 
Let $\Sigma^{(1)}=\mathbb{C}/{\Lambda}$ denote an oriented torus 
for lattice ${\Lambda}=2\pi i(\Z\tau\oplus \Z)$ and $\tau\in \mathbb{H}_{1}$.    
Define annuli $\mathcal{A}_{a}$, $a=1$, $2$ centered at $z=0$ and $z=w$ of  $\Sigma^{(1)}$ with local
coordinates $z_{1}=z$ and $z_{2}=z-w$ respectively.
Define the following label convention 
\begin{equation}
\overline 1=2,\quad \overline 2=1.
\label{bardefn}
\end{equation}   
Take the outer radius of 
$\mathcal{A}_{a}$ to be $r_{a}<\half D(q)$ for $D(q)=\min_{\lambda \in {\Lambda}, \lambda \neq 0}|\lambda |$ 
and the inner radius
to be $|{\rho }|/r_{\overline{a}}$, with $|\rho|\leq r_{1}r_{2}$, 
where 
we introduce a complex parameter $\rho$, $|\rho|\le r_1 r_2$. 
$r_{1}$, $r_{2}$ must be sufficiently small to ensure that the disks do not
intersect.
 Excise the disks
\begin{equation*}
\{z_{a}, \left\vert z_{a}\right\vert <|\rho |r_{\overline{a}}^{-1}\}\subset 
\Sigma^{(1)}, 
\end{equation*} 
to form a twice-punctured surface 
\begin{equation*}
\widehat{\Sigma}^{(1)}=\Sigma^{(1)}\backslash \bigcup_{a=1,2}\{z_{a},  \left\vert
z_{a}\right\vert <|\rho |r_{\overline{a}}^{-1}\}.
\end{equation*}
 We define annular regions 
$\mathcal{A}_{a}\subset \widehat{\Sigma}^{(1)}$ with $\mathcal{A}_{a}=
\{z_{a}, |{\rho }|r_{\overline{a}}^{-1}\leq \left\vert z_{a}\right\vert \leq r_{a}\}$ and
identify them as a single region $\mathcal{A}=\mathcal{A}_{1}\simeq \mathcal{A}_{2}$ 
via the sewing relation 
\begin{equation}
 \label{rhosew}
z_{1}z_{2}=\rho,
\end{equation}
to form a compact genus two Riemann surface 
$\Sigma^{(2)}=\widehat{\Sigma}^{(1)}\backslash \{\mathcal{A}_{1}\cup \mathcal{A}_{2}\}\cup \mathcal{A}$, 
parameterized by \cite{MT2}
\begin{equation}
\mathcal{D}^{\rho }=\{(\tau ,w,\rho )\in \mathbb{H}_{1}\times \mathbb{C}
\times \mathbb{C}\ , |w-\lambda |>2|\rho |^{\half }>0,\ \lambda \in {\Lambda}\}.
\label{Drho}
\end{equation}
\subsection{The Genus Two Szeg\H{o} Kernel in the $\rho$-Formalism} 
\label{relationslambda1}
Now we recall the construction of the Szeg\H{o} kernel on a genus two Riemann surface 
in the $\rho$-formalism \cite{TZ1}. 
The genus one prime form for $x,y\in \mathbb{C}$ and $\tau \in \mathbb{H}_{1}$ is
\begin{eqnarray}
E(x,y)&=&K (x-y,\tau )\; dx^{-\half }\; dy^{-\half },
\label{Prime1}
\end{eqnarray}
where 
\begin{equation}
K (z,\tau )=\frac{\vartheta _{1}(z,\tau )}{\partial_{z}\vartheta_{1}(0,\tau)},
\quad  \vartheta_{1}(z,\tau) 
= \vartheta 
\begin{bmatrix}
\half  
\\
\half 
\end{bmatrix}
(z,\tau ). 
\label{Ktheta}
\end{equation}
Let $(\theta_1 ,\phi_1 )\in U(1)\times U(1)$ with $(\theta_1 ,\phi_1 )\neq (1,1)$.
The genus one  Szeg\H{o} kernel is
\begin{equation}
S^{(1)}
\begin{bmatrix}
\theta_1  \\ 
\phi_1 
\end{bmatrix}
 (x,y\; \vert \tau )=
P_{1}
\begin{bmatrix}
\theta_1  \\ 
\phi_1 
\end{bmatrix}
(x-y,\tau )\; dx^{\half }\; dy^{\half },
\label{S1}
\end{equation}%
where
\begin{eqnarray}
P_{1}
\begin{bmatrix}
\theta_1  \\ 
\phi_1 
\end{bmatrix} (z,\tau )&=&
\frac{\vartheta 
\begin{bmatrix}
\alpha_1 \\ 
\beta_1
\end{bmatrix} (z,\tau )}
{\vartheta 
\begin{bmatrix}
\alpha_1 
\\ 
\beta_1
\end{bmatrix}(0,\tau )}\frac{1}{K (z,\tau)},
\label{p_one_theta}
\end{eqnarray}%
where $-\phi_1=\exp (2\pi i \alpha_1)$ and $-\theta_1 =\exp (-2\pi i \beta_1)$   are the periodicities of $S^{(1)}
\left[
\begin{smallmatrix}
\theta_1  \\ 
\phi_1 
\end{smallmatrix}
\right]
 (x,y \vert  \tau )$ in $x$ on the standard $a$ and $b$ cycles respectively.

In  \cite{TZ1} we determine the genus two Szeg\H{o} kernel
\[
S^{(2)}(x,y) = S^{(2)}
 \begin{bmatrix}
{\theta}^{(2)} \\
 {\phi}^{(2)}
\end{bmatrix} 
 (x,y),
 \]
with periodicities $\left({\theta}^{(2)}, {\phi}^{(2)} \right)= ({\theta}_{i}, {\phi}_{i})$ for $i=1,2$ 
on a natural homology basis on the genus two Riemann surface $\Sigma^{(2)}$ 
formed by self-sewing the torus $\Sigma^{(1)}$ in terms of genus one Szeg\H{o} kernel data 
$S^{(1)}(x,y)=S^{(1)}\big[ \begin{smallmatrix}{\theta_1} \\{\phi_1} \end{smallmatrix} \big](x,y)$ with
$(\theta_1,\phi_1) \neq (1,1)$.  
 The $S^{(2)}$ multipliers \eqref{eq:periodicities} 
on the cycles $a_1$, $b_1$ are determined by the corresponding multipliers of $S^{(1)}$ so that
$\phi_1^{(2)}=\phi_1$ and $\theta_1^{(2)}=\theta_1$ 
i.e. $\alpha_1^{(2)}=\alpha_1$ and $\beta_1^{(2)}=\beta_1$.  
The remaining  multipliers 
$\phi_{2}=\phi_{2}^{(2)}=-e^{2\pi i\alpha_{2}^{(2)}}$ and
$\theta_{2}=\theta_{2}^{(2)}=-e^{-2\pi i\beta_{2}^{(2)}}$ 
on the cycles $a_{2}$ and $b_{2}$ must be specified by the additional conditions
\begin{eqnarray*}
 S^{(2)}(e^{2\pi i}x_a,y) &=& -\phi_{2}^{ a-\overline{a}} \; 
 S^{(2)}(x_{a},y),    \label{agplus_mult}
\\
 S^{(2)}(x_a,y) &=& -\theta_{2}^{ a-\overline{a}} \; 
 S^{(2)}(x_{\overline{ a}},y),    \label{bgplus_mult}
\end{eqnarray*}
for $x_{a}\in \mathcal{A}_{a}$ and $x_{\overline{a}}\in\mathcal{A}_{\overline a}$ (recalling the label convention \eqref{bardefn}).

It is convenient to define $\kappa\in \left(\half ,\half \right)$ by
$\phi_{2}=-e^{2\pi i\kappa}$ i.e. $\kappa=\alpha_{2}^{(2)} \mod 1$. 
We then find that 
 $S^{(2)}$ is holomorphic in $\rho $ for $|\rho|<r_{1}r_{2}$ with \cite{TZ1} 
\begin{equation*}
S^{(2)}(x,y)=S_{\kappa}(x,y)+O(\rho),
\label{eq:SgpSkappa}
\end{equation*} 
for $x$, $y\in \widehat{\Sigma}^{(1)}$ and
\begin{eqnarray}
\label{skappa_1}
&&S_{\kappa}(x,y)=
\left(
\frac {\vartheta_{1} (x-w,\tau) \vartheta_{1} (y,\tau)} 
{ \vartheta_{1} (x,\tau)  \vartheta_{1} (y-w,\tau)}
\right)^{\kappa}\,
\frac{ \vartheta  \begin{bmatrix}{\alpha_{1} } \\ {\beta_{1} }\end{bmatrix} 
\left( x-y +\kappa w,\tau\right)}
{\vartheta  \begin{bmatrix}{\alpha_{1} } \\ {\beta_{1} }\end{bmatrix}  
\left(\kappa w,\tau \right)  K (x-y,\tau)} 
 dx^{\half} dy^{\half}, \qquad
\end{eqnarray}   
for $\kappa\neq -\half $
(with a different expression  when $\kappa=- \half $ given in \cite{TZ1}). 
We will assume  $\kappa\neq-\half$ throughout this paper. 
Note also  that $S_{\kappa=0}(x,y)=S^{(1)}
\left[
\begin{smallmatrix}
\theta_1  \\ 
\phi_1 
\end{smallmatrix}
\right]
 (x,y)$, the genus one Szeg\H{o} kernel.

$S_\kappa(x,y)$ has an  expansion in the neighborhood of the punctures at $0,w$ in terms of local coordinates $x_1=x,y_1=y$ and $x_1=x-w,y_2=y-w$ as follows \cite{TZ3}
\begin{eqnarray}
S_{\kappa}(x_{\abar},y_b)
=
\left[
\delta_{\abar,b}
\frac{1}{x_{b}-y_b}
\left(\frac{x_b}{y_b}\right)^{\kappa(-1)^{b}}
+\sum\limits_{k,l\ge 1}C_{a b}(k,l)x_{\abar}^{k_a-1}y_b^{l_b-1}
\right]
dx_{\abar}^{\half } dy_b^{\half },\quad
\label{Skexp}
\end{eqnarray}
where $C_{a b}(k,l)= 
C_{a b}\begin{bmatrix}{\theta_{1} } \\ {\phi_{1} }\end{bmatrix}(k,l\vert \tau,w,\kappa)$ and  $k_a =k  +\kappa(-1)^{\overline{a}}$  for integer $k\ge 1$ and $a=1,2$.
We may invert this to obtain the infinite block moment matrix 
\begin{eqnarray}
&& C_{a b}(k,l)= \frac{ 1} { (2\pi i)^2 }  
\oint_{\mathcal{C}_{{\overline{a}}}(x_{\overline{a}})}
\oint_{\mathcal{C}_{b}(y_b)}
(x_{\overline{a}})^{-k_a} (y_b)^{-l_b}
S_{\kappa} (x_{\overline{a}},y_b) 
\; dx_{\overline{a}}^{\half } \; dy_b^{\half }.\quad
\label{Cijdef}
\end{eqnarray}
We also define half-order differentials 
\begin{eqnarray}
d_{a}(x,k)
&=&  
\frac{1}{2\pi i}
\oint_{\mathcal{C}_{a}(y_{a})}
y_{a}^{-k_a} 
S_{\kappa}(x, y_{a})  
dy_{a}^{\half },  
\label{ddef1}\\
\overline{d}_{a}(y,k)
&=&
\frac{1}{2\pi i}
\oint_{\mathcal{C}_{\overline{a}}(x_{\overline{a}})}
x_{\overline{a}}^{-k_a} 
S_{\kappa}(x_{\overline{a}}, y)  
dx_{\overline{a}}^{\half }.
\label{ddef2}
\end{eqnarray}
We further define
\begin{eqnarray}
G_{a b}(k,l)&=&  \rho^{\half (k_a+l_b-1) }C_{a b}(k,l),  
\label{Gijdef}\\
h_{a}(x,k)&=&\rho^{\half (k_{a} - \half )}d_{a}(x,k),
\label{hdef1}
\\
\overline{h}_{a}(y,k)&=&\rho^{\half (k_{a} - \half )}\overline{d}_{a}(y,k),
\label{hdef2}
\end{eqnarray} 
with associated infinite matrix $G=\left(G_{ab}(k,l)\right)$ and row vectors
$h(x)=(h_{a}(x,k))$ and $\overline{h}(y)=(\overline{h}_{a}(y,k))$. 
From the sewing relation \eqref{rhosew} we have
\begin{equation}
dz_a^{\half }=(-1)^{\overline a}\; \xi\; \rho^{\half } \; \frac{dz_{\overline a}^{\half }}{z_{\overline a}}, 
\label{dz1dz2rho}
\end{equation}
for $\xi\in\{\pm \sqrt{-1}\}$ 
  depending on the branch of the double cover of $\Sigma^{(1)}$ chosen. 
Finally define 
\begin{equation}
\label{T_matrix}
T=\xi GD^{\theta_2}, 
\end{equation}
with an infinite diagonal matrix 
\begin{equation}
\label{Dtheta}
D^{\theta_2}(k,l)=
\left[
\begin{array}{cc}
\theta_2^{-1} & 0
\\
0 & -\theta_2 
\\
\end{array}
\right]
\delta(k,l).
\end{equation}
Defining $\det (I - T)$ by 
$\log\det (I - T) = \Tr \log\left(I - T\right)=-\sum_{n\ge 1}\frac{1}{n} \Tr T^n$, we find
\begin{theorem}[\cite{TZ1}]
\label{theorem_S2}
On the $\rho$-sewing domain $\mathcal{D}^\rho$ we have:
\begin{enumerate}
[(i)]
\item
$(I-T)^{-1}=\sum_{n\ge 0} T^n$ is convergent.
\item 
$\det (I - T)$
is holomorphic and non-vanishing.
\item
The genus two Szeg\H{o} kernel for $x,y\in\widehat{\Sigma}^{(1)}$ is given by
\begin{equation}
\label{s_genus_two}
S^{(2)}(x,y) =S_{\kappa}(x,y) +\xi h(x)D^{\theta} (I-T)^{-1} \overline{h}^{T}(y),      
\end{equation}
where $\overline{h}^{T}(y)$ denotes the transpose of the infinite row vector $\overline{h}(y)$.
\end{enumerate}
\end{theorem}

\section{The Free Fermion VOSA and its Twisted Modules} 
\label{vosas}
\subsection{Vertex Operator Superalgebras}
\label{vos}
We review some aspects of Vertex Operator Superalgebra  theory (e.g. 
\cite{FHL, FLM, K, MN, MT5}).
We define the standard formal series 
\begin{eqnarray}
\delta\left(\frac{x}{y}\right)&=&\sum_{r\in \Z} x^{r}y^{-r},
\label{delta}\\
(x+y)^{\lambda}&=&\sum_{s\ge 0}\binom{\lambda}{s}x^{\lambda-s}y^s,
\label{xylambda}
\end{eqnarray} 
for any formal variables $x,y,\lambda$ where $\binom{\lambda}{s}=\frac{\lambda(\lambda-1)\ldots (\lambda-s+1)}{s!}$.  
\begin{definition}\label{defvosa}
A Vertex Operator Superalgebra (VOSA) is determined by a quadruple 
$(V,Y,\mathbf{1},\omega )$ 
as follows: $V$ is a superspace $V=V_{\bar{0}}\oplus V_{\bar{1}}$
with parity $p(u)=0$ or $1$ for $u\in V_{\bar{0}}$ or $V_{\bar{1}}$ respectively. $V$ also has a $\half\Z$-grading with $V=\bigoplus _{r \in \half\Z}V_{r}$  
with $\dim V_{r}<\infty$.
   $\vac\in V_{0}$ is the vacuum vector and 
$\omega\in V_{2}$ is the conformal vector with properties described below.
 
$Y$ is a linear map $Y:V\rightarrow \mathrm{End}(V)[[z,z^{-1}]]$ for formal
variable $z$ so that for any vector $u\in V$ we have a vertex operator  
\begin{equation}
Y(u,z)=\sum_{n\in \Z}u(n)z^{-n-1}.  
\label{Ydefn}
\end{equation}
The linear operators (modes) $u(n):V\rightarrow V$ satisfy creativity 
\begin{equation}
Y(u,z)\vac = u +O(z),
\label{create}
\end{equation}
and lower truncation 
\begin{equation}
u(n)v=0,
\label{lowertrun}
\end{equation}
for each $u,v\in V$ and $ n\gg 0$.
For the conformal vector $\omega$ 
\begin{equation}
Y(\omega ,z)=\sum_{n\in \Z}L(n)z^{-n-2},  \label{Yomega}
\end{equation}
where $L(n)$ satisfies the Virasoro algebra for some central charge $c$ 
\begin{equation}
[\,  L(m),L(n)\, ]=(m-n)L(m+n)+\frac{c}{12}(m^{3}-m)\delta_{m,-n}\Id_V.
\label{Virasoro}
\end{equation}
Each vertex operator satisfies the translation property 
\begin{equation}
Y(L(-1)u,z)=\delz Y(u,z).  
\label{YL(-1)}
\end{equation}
The Virasoro operator $L(0)$ provides the $\half\Z$-grading with $L(0)u=ru$ for 
$u\in V_{r}$ and with $r\in \Z+\half p(u)$. 
Finally, the vertex operators satisfy the Jacobi identity
\begin{eqnarray}
\notag
& z_0^{-1}\delta\left( \frac{z_1 - z_2}{z_0}\right) Y (u, z_1 )Y(v , z_2)    
  -(-1)^{p(u)p(v)}z_0^{-1} \delta\left( \frac{z_2 - z_1}{-z_0}\right) Y(v, z_2) Y(u , z_1 ) 
&\\
\label{VOAJac}
& 
= z_2^{-1}    
\delta\left( \frac{z_1 - z_0}{z_2}\right)
Y \left( Y(u, z_0)v, z_2\right).  &
\end{eqnarray} 
$(V,Y,\mathbf{1},\omega )$ is called a Vertex Operator Algebra (VOA) when $V_{\bar{1}}=0$.
\end{definition}

These axioms imply 
locality, skew-symmetry, associativity and commutativity:
\begin{eqnarray}
(z_{1}-z_{2})^N
Y(u,z_{1})Y(v,z_{2}) 
&=& (-1)^{p(u)p(v)}(z_{1}-z_{2})^N
Y(v,z_{2})Y(u,z_{1}),
\notag
\\
\label{Local}\\
Y(u,z)v &=& (-1)^{p(u)p(v)} e^{zL(-1)}Y(v,-z)u,
\label{skew}
\\
(z_{0}+z_{2})^N Y(u,z_{0}+z_{2})Y(v,z_{2})w &=& (z_{0}+z_{2})^N Y(Y(u,z_{0})v,z_{2})w,
\label{Assoc}\\
u(k)Y(v,z)-(-1)^{p(u)p(v)} Y(v,z)u(k)
&=& \sum_{j\ge 0}\binom{k}{j}
Y(u(j)v,z)z^{k-j},\label{Comm}
\end{eqnarray}
for $u,v,w\in V$ and integers $N\gg 0$  \cite{FHL, K, MT5}.

We next define the notion of a twisted $V$-module. 
Let $g$ be a $V$-automorphism $g$ i.e.~a linear map preserving $\vac$ and $\omega$ such that 
\begin{equation*}
g Y(v,z) g^{-1}=  Y(g v,z),
\end{equation*}
for all $v \in V$. We assume that $V$ can be decomposed into $g$-eigenspaces
\begin{equation*}
V= \oplus_{\rho\in \mathbb{C}} V^{\rho}, 
\end{equation*}
where $V^{\rho}$ denotes the eigenspace of $g$ with eigenvalue $e^{2\pi i \rho}$.
\begin{definition}\label{Vgmodule}
A $g$-twisted $V$-module for a VOSA $V$ is a pair $(W^{g}, Y_{g})$ 
where $W^{g}$ is a $\mathbb{C}$-graded vector space
$W^{g}=\bigoplus\limits_{r\in \mathbb{C}}W^{g}_{r}$ with $\dim W_{r}<\infty$ and where    $W_{r+n}=0$ for all $r$ and  $n\ll 0$. 
$Y_{g}$ is a linear map $Y_{g}:V\rightarrow \mathrm{End\ }W^{g}\{z\}$, the vector space of 
$\mathrm{End\ }W^{g}$-valued formal series in $z$ with arbitrary complex powers of $z$. Then for $v\in V^{\rho}$
\begin{equation*}
Y_{g}(v,z)= \sum_{n \in \rho+ \mathbb {Z}} v_{g}(n)z^{-n-1}, 
\end{equation*}
with $v_{g}(\rho+l)  w=0$ for $w\in W^{g}$ and $l \in \mathbb {Z}$ sufficiently large. 
For the vacuum vector $Y_{g}(\vac, z) = \mathrm{Id}_{W^{g}}$ and for the conformal vector
\begin{equation}
Y_{g}(\omega ,z)=\sum_{n\in \Z}L_{g}(n)z^{-n-2},  \label{Ygomega}
\end{equation} 
where $L_{g}(0)w=r w$ for $w\in W^{g}$. The $g$-twisted vertex operators satisfy the twisted Jacobi identity: 
\begin{eqnarray}
\nonumber
&   z_0^{-1} \delta\left( \frac{z_1 - z_2}{z_0}\right) Y_{g}(u, z_1) Y_{g}(v, z_2)  
- (-1)^{p(u,v)} z_0^{-1} \delta\left( \frac{z_2 - z_1}{-z_0}\right) Y_{g}(v, z_2) Y_{g}(u, z_1)  
&\\
\label{gJac}
& \qquad =
 z_2^{-1} \left( \frac{z_1 - z_0}{-z_2}\right)^{-\rho} 
\delta\left( \frac{z_1 - z_0}{-z_2}\right) Y_{g}(Y(u, z_0)v, z_2), &
\end{eqnarray}
for $u \in V^{\rho}$.  
\end{definition}
This definition is an extension of the standard one for $g$ unitary where $\rho\in \mathbb{R}$ \cite{DLinM}. These axioms  imply that $L_{g}(n)$ of \eqref{Ygomega} satisfies the Virasoro algebra \eqref{Virasoro} for the same central charge $c$ and that the translation property holds: 
\begin{eqnarray}
Y_{g}(L(-1)u,z)&=&\delz Y_{g}(u,z).
\label{YL(-1)g}
\end{eqnarray} 
\subsection{The Free Fermion VOSA}
\label{FVOSA}
We consider in this paper the rank two free fermion VOSA 
$V(H,\Z+\half )^{\otimes 2}$ of central charge 1. 
The weight $\half$ space $V_{\half}$ is spanned by $\psi^{+},\psi^{-}$ with vertex operators
\[
Y(\psi^{\pm },z)=\sum_{n\in \Z }
\psi ^{\pm }(n)z^{-n-1},
\]
whose modes satisfy the anti-commutation relations 
 \begin{equation}
 \lbrack \psi ^{+}(m),\psi ^{-}(n)\rbrack=
 \delta _{m,-n-1}, \quad \lbrack \psi^{+}(m),\psi^{+}(n)\rbrack
 =\lbrack \psi^{-}(m),\psi^{-}(n)\rbrack=0. 
\label{psiplus_minus_comm}
 \end{equation}
The VOSA is generated by $\psi^{\pm}$ with
$V$ spanned by Fock vectors of the form
\begin{equation}
\label{Fock}
\Psi(\mathbf{k},\mathbf{l})\equiv \psi^{+}(-k_{1})\ldots \psi^{+}
(-k_{s}) 
\; \psi^{-}(-l_{1})
 \ldots \psi^{-}(-l_{t}) \vac, 
\end{equation}
for distinct $0<k_{1}<\ldots <k_{s}$ and $0<l_{1}<\ldots <l_{t}$.
The Virasoro vector
\[
\omega=\half(\psi^{+}(-2)\psi^{-}(-1)+\psi^{-}(-2)\psi^{+}(-1))\vac.
\]
generates a Virasoro algebra with central charge $c=1$ for which the Fock vectors have weight
\begin{eqnarray}
\wt(\Psi(\mathbf{k},\mathbf{l}))=  \sum_{i=1}^s 
\left(k_i-\half\right) + \sum_{j=1}^t \left(l_j -\half \right) .
\label{eq:Psswt}
\end{eqnarray}


\subsection{Bosonization}
The weight $1$ space $V_1$ is spanned by the Heisenberg vector
$a=\psi^{+}(-1)\psi^{-}$ whose modes obey the Heisenberg commutation relations 
\[
[a(m),a(n)]=m\delta_{m,-n}.
\] 
Then $\omega=\half a(-1)^{2}\vac$ is the usual conformal vector for the Heisenberg
VOA $M$ with unit central charge.   
As is well known, we may decompose $V$ into irreducible $M$-modules $M\otimes e^{m}$ 
with $a(0)$ eigenvalue $m\in \Z$ so that $V(H,\Z+\half )^{\otimes 2}\cong V_{\Z}=\oplus_{m\in \Z} M\otimes e^{m}$, the lattice VOSA for the $\Z$-lattice with trivial cocycle structure. 
In particular, the generators $\psi^{+},\psi^{-}$ are given by 
\begin{eqnarray}
\psi^{+}=\vac\otimes e^{1},\quad \psi^{-}=\vac\otimes e^{-1}.
\label{eq: psipmlat}
\end{eqnarray}

In terms of Heisenberg modes, the vertex operators for $u\otimes e^m\in V_{\Z}$ for $u\in M$ are given by (e.g. \cite{K}, \cite{TZ3})
\begin{eqnarray}
Y(u\otimes e^{m},z)&=& e^{m \widehat{q}}\, Y_{-}(m, z)  
Y(u ,z)Y_{+}(m,z)\,z^{ma(0)},
\label{Ym}
\end{eqnarray}
for $m\in \Z$ where
\begin{equation}
\label{pcom1}
\left[\,a(n), \widehat{q} \,\right] =  \delta_{n,0}, 
\end{equation}  
and
\begin{eqnarray}
\label{Ypm}
Y_{\pm}(m, z) &=& \exp \left(\mp \,m  \sum_{n>0}\frac{ a(\pm \; n)}{n} z^{\mp n}\right). 
\end{eqnarray}

\subsection{$g$-Twisted $V_{\Z}$-Modules and a Generalized VOA}
The mode $a(0)$ generates continuous $V_{\Z}$-automorphisms 
$g=e^{-2\pi i \alpha a(0)}$ for all $\alpha\in \mathbb{C}$. In particular, we define the fermion number involution
\[
\sigma=e^{\pi i a(0)},
\]
where $\sigma u=(-1)^{p(u)}u$ for $u$ of parity $p(u)$.
Define for all $\alpha\in \mathbb{C}$
\begin{eqnarray}
Y_{g}(u,z)&=Y(\Delta(\alpha,z)u,z),\label{Y_g}\\
\Delta(\alpha,z)&=z^{\alpha a(0)} Y_{+}\left (\alpha, -z\right), \label{Delta}
\end{eqnarray}
where
\begin{eqnarray}
\label{Ypmalpha}
Y_{\pm}(\alpha, z) &=& \exp \left(\mp \,\alpha  \sum_{n>0}\frac{ a(\pm \; n)}{n} z^{\mp n}\right). 
\end{eqnarray}
We note the identity 
\begin{eqnarray}
Y_{+}(\alpha,x)Y_{-}(\beta,y) &=& \left(1-\frac{y}{x}  \right)^{\alpha\beta}Y_{-}(\beta,y) Y_{+}(\alpha,x).
\label{Ypmcom}
\end{eqnarray}
Then we have (see \cite{Li} for details)
\begin{proposition}
$(V_{\Z},Y_{g})$ is a $g$-twisted $V_{{\Z}}$-module.
\end{proposition} 
Thus the twisted module is constructed by the action of the twisted vertex operators $Y_g(u,z)$ on the original vector space $V_{\Z}$.
In Section~5 of \cite{TZ3} an isomorphic construction is described whereby the $g$-twisted module is determined by the action of the original vertex operators \eqref{Ym} on a twisted vector space $ V_{{\Z}+\alpha}=e^{\alpha \widehat{q}}V_{\Z}=\oplus_{m\in \Z} M\otimes e^{m+\alpha}$ where 
\begin{eqnarray}
Y_g(u,z)=e^{-\alpha \widehat{q}}Y(u,z)e^{\alpha \widehat{q}},
\label{Yg2}
\end{eqnarray}
for all $u\in V_{\Z}$. In particular, for the Virasoro vector, we find $Y_g(\omega,z)$ has $g$-twisted modes
\begin{eqnarray}
L_{g}(n)=L(n)+\alpha a(n)+\half \alpha^2 \delta_{n,0}.
\label{eq:Lgn}
\end{eqnarray}
Thus the $g$-twisted grading operator is $L_{g}(0)=L(0)+\alpha a(0)+\half \alpha^2$ so that $u\otimes e^m\in V_{\Z}$ has
 $g$-twisted Virasoro weight $\wt(u)+\half(m+\alpha)^2$ equal to the $L(0)$ Virasoro weight of  $u\otimes e^{m+\alpha}\in V_{\Z+\alpha}$. 
Hence, the Fock vector $\Psi(\mathbf{k},\mathbf{l})\in M\otimes e^{s-t}$ has $g$-twisted Virasoro weight $\wt(\Psi(\mathbf{k},\mathbf{l}))+\alpha (s-t)+\half \alpha^2$ equal to the $L(0)$ weight of the $V_{\Z+\alpha}$ twisted Fock vector
\begin{eqnarray}
\Psi_{\alpha}(\mathbf{k},\mathbf{l})\equiv e^{\alpha \widehat{q}}\Psi(\mathbf{k},\mathbf{l})\in M\otimes e^{s-t+\alpha}.
\label{eq:Focka}
\end{eqnarray}

In \cite{TZ3} we describe a Generalized VOA formed by the Heisenberg VOA $M$ extended by vertex operators (creative intertwiners) \emph{for all} the irreducible $M$-modules $M\otimes e^{\alpha}$. Choosing a trivial cocycle structure, as per the standard VOSA $V_{\Z}$, we define for all $\alpha\in \mathbb{C}$
\begin{eqnarray}
{\cal Y}(u\otimes e^{\alpha},z)&=& e^{\alpha \widehat{q}}\, Y_{-}(\alpha, z)  
Y(u ,z)Y_{+}(\alpha,z)\,z^{\alpha a(0)}.
\label{Yalpha}
\end{eqnarray}
Clearly, ${\cal Y}(u\otimes e^{m},z)=Y(u\otimes e^{m},z)$ for all $m\in \Z$.  We also note \cite{TZ3} 
\begin{eqnarray}
e^{-\alpha \qhat}\Ycal(u\otimes e^{\beta},z)e^{\alpha\qhat}=\Ycal(\Delta(\alpha,z) u\otimes e^{\beta},z),
\label{eq:qhatconj}
\end{eqnarray}
which implies \eqref{Yg2}. 

The vertex operators \eqref{Yalpha} possess the creativity and translation property with Heisenberg Virasoro vector and satisfy a generalized Jacobi identity (with a trivial cocycle structure) as follows:
\begin{theorem}[\cite{TZ3}]
\label{GenVOA} 
 $(\Mcal,\Ycal,\vac,\omega)$ is a Generalized VOA with vertex operators ${\cal Y} (u\otimes e^{\alpha}, z )  \in\Lin({\cal M},{\cal M})[[z,z^{-1}]]$ obeying the generalized Jacobi identity  
\begin{eqnarray}
\notag
& & z_0^{-1} \left( \frac{z_1 - z_2}{z_0}\right)^{ 
-\alpha \beta
}  
\delta\left( \frac{z_1 - z_2}{z_0}\right)  
\; {\cal Y}  \left(u\otimes e^{\alpha}, z_1 \right) \; 
{\cal Y}  \left(v \otimes e^{\beta}, z_2\right)   
\\
\notag
& & \; 
  - z_0^{-1}
\left( \frac{z_2 - z_1}{z_0}\right)^{ -\alpha\beta 
}  
 \delta\left( \frac{z_2 - z_1}{-z_0}\right)
 \;{\cal Y} \left(v \otimes e^{\beta}, z_2\right)  \; {\cal Y} \left(u\otimes e^{\alpha}, z_1 \right) 
\\
\label{Jac}
&& 
\;
= z_2^{-1}  
 \;   
\delta\left( \frac{z_1 - z_0}{z_2}\right)
{\cal Y} \left( {\cal Y} \left(u \otimes e^{\alpha}, z_0\right)
 (v\otimes e^{\beta}), z_2\right)
\left( \frac{z_1 - z_0}{z_2}\right)^{\alpha a(0)},\quad 
\end{eqnarray} 
for all $u\otimes e^{\alpha},v\otimes e^{\beta}\in{\cal M}$.
\end{theorem}
We exploit the following features of this Generalized VOA: (i) generalized locality,  skew-symmetry, associativity and Heisenberg commutivity and  (ii) an invariant bilinear form on ${\cal M}$ which we now describe.
\subsubsection{Generalized Locality,  Skew-Symmetry, Associativity and Heisenberg Commutivity}
The Generalized VOA of \eqref{Jac} leads to a more general notion of locality, skew-symmetry and  associativity  \eqref{Local}--\eqref{Assoc}. 
We also note a  generalization of \eqref{Comm} where commutation with respect Heisenberg modes are only considered.
\begin{proposition}
\label{prop_assoc}
For $u\otimes e^{\alpha},v\otimes e^{\beta},w\otimes e^{\gamma} \in{\cal M}$ and 
for integer $N\gg 0$
\begin{eqnarray}
& &
\notag
\left(z_1-z_2\right)^{N}\left(z_1-z_2\right)^{-\alpha \beta}
\; {\cal Y}  \left(u\otimes e^{\alpha}, z_1 \right) 
{\cal Y}  \left(v \otimes e^{\beta}, z_2\right) 
\\
& &
=
\left(z_1-z_2\right)^{N}\left(z_2-z_1\right)^{-\alpha \beta}\; 
{\cal Y}  \left(v \otimes e^{\beta}, z_2\right)
{\cal Y}  \left(u\otimes e^{\alpha}, z_1 \right),\;
\label{genloc}
\\
& &
z^{-\alpha \beta}
\; {\cal Y}  \left(u\otimes e^{\alpha}, z\right) v \otimes e^{\beta}
=
(-z)^{-\alpha \beta}\; e^{zL(-1)}\;{\cal Y}  \left(v\otimes e^{\beta}, -z\right) u \otimes e^{\alpha},\quad
\label{genskew}
\\
& &
\notag
\left(z_0+z_2\right)^{N-\alpha \gamma}
\; {\cal Y}  \left(u\otimes e^{\alpha}, z_0+z_2 \right) 
{\cal Y}  \left(v \otimes e^{\beta}, z_2\right) \,  w\otimes e^{\gamma}
\\
& &=
\left(z_2+z_0\right)^{N-\alpha \gamma} \;
{\cal Y} \left( {\cal Y} \left(u \otimes e^{\alpha}, z_0\right)
 (v\otimes e^{\beta}), z_2\right)
\,  w\otimes e^{\gamma},\;
\label{genassoc}
\\
&&
\left[u(k),\Ycal( v \otimes e^{\beta},z)\right]
= \sum_{j\ge 0}\binom{k}{j}
\Ycal(u(j) v \otimes e^{\beta},z)z^{k-j},\quad u\in M.
\label{genComm}
\end{eqnarray} 
\end{proposition} 
The proof appears in Appendix \ref{PropGen}.

\medskip
It is convenient to define for formal parameter $z$ and  $\chi\in \C$
\begin{equation}
(-z)^{\chi}=e^{i\pi B\chi } z^{\chi},
\label{minz}
\end{equation}
where $B$ is an \emph{odd integer} parametrizing the formal branch cut\footnote{Note some notational changes from \cite{TZ3}}. Then generalized locality and skew-symmetry can be rewritten as 
\begin{eqnarray}
& &
\notag
\left(z_1-z_2\right)^{N-\alpha \beta}
\; {\cal Y}  \left(u\otimes e^{\alpha}, z_1 \right) 
{\cal Y}  \left(v \otimes e^{\beta}, z_2\right) 
\\
& &
=
e^{-i\pi B\alpha\beta}\left(z_1-z_2\right)^{N-\alpha \beta}\; 
{\cal Y}  \left(v \otimes e^{\beta}, z_2\right)
{\cal Y}  \left(u\otimes e^{\alpha}, z_1 \right),\;
\label{genlocB}
\\
& &
{\cal Y}  \left(u\otimes e^{\alpha}, z\right) v \otimes e^{\beta}
=
e^{-i\pi B\alpha\beta}\; e^{zL(-1)}\;{\cal Y}  \left(v\otimes e^{\beta}, -z\right) u \otimes e^{\alpha}.\quad
\label{genskewB}
\end{eqnarray} 
These reduce to \eqref{Local} and \eqref{skew} for $\alpha,\beta\in\Z$.
\subsubsection{The Invariant Form on ${\cal M}$}
\label{liza}
In \cite{TZ3} we introduced an invariant bilinear form $\langle \cdot , \cdot\rangle$ on $\Mcal$ 
associated with the M\"obius map 
\cite{FHL}, \cite{S}, \cite{TZ2}
\begin{equation}
\left(
\begin{array}{cc}
0 & \lambda\\
-e^{ i\pi B } \lambda^{-1} & 0\\	
\end{array}
\right)
:z\mapsto -\frac{\lambda^{2}}{e^{ i\pi B } z},
 \label{eq: gam_lam}
\end{equation}
for $\lambda\neq 0$. 
We are particularly interested in the M\"{o}bius
 map $z\mapsto \rho/z$ associated with the sewing condition \eqref{rhosew} so that we will later choose
\begin{equation}
\lambda=e^{\half i\pi B } \rho^{\half },
\label{eq:lamb_eps}
\end{equation}
for the odd integer $B$ of \eqref{minz}. 
Thus we reformulate the sewing relationship \eqref{rhosew} as $z_1=-\frac{\lambda^{2}}{ z_2}$ so that \eqref{dz1dz2rho} implies $dz_1^{\half }= \xi \rho^{\half }/z_2\, dz_{2}^{\half }$ for $\xi=e^{\half i\pi B }$ .

Define the adjoint vertex operator
\begin{eqnarray}
{\cal Y}^\dagger\left(u\otimes e^{\alpha}, z\right) &=&
{\cal Y}\left(
e^{-z\lambda^{-2}L(1)} 
\left(\frac{\lambda}{e^{i\pi B} z}\right)^{2 L(0)}
(u\otimes e^{\alpha}),\frac{\lambda^{2}}{e^{i\pi B} z}\right). 
\label{eq: calY_dag}
\end{eqnarray}
A bilinear form $\langle \cdot , \cdot\rangle$ on ${\cal M}$ is said to be   
invariant if for all $u\otimes e^{\alpha}$, $v\otimes e^{\beta}$, $w\otimes e^{\gamma} \in {\cal M}$ we have
\begin{equation}
\langle {\cal Y}(u\otimes e^{\alpha},z) \; (v\otimes e^{\beta}), \; w\otimes e^{\gamma}\rangle 
= 
e^{-i\pi B \alpha\beta}  
\langle v\otimes e^{\beta}, \; {\cal Y}^{\dagger}(u\otimes e^{\alpha},z)\; w\otimes e^{\gamma}\rangle. 
\label{eq: inv bil form1}
\end{equation}
\eqref{eq: calY_dag} reduces to the usual definition for a VOSA when $\alpha,\beta,\gamma\in\Z$  \cite{S, TZ2}. 
Choosing the normalization $\langle {\bf 1},{\bf 1} \rangle=1$  then $\langle \cdot , \cdot \rangle$ on $\cal{M}$ is symmetric, unique and invertible with \cite{TZ3}
\begin{equation}
\langle u\otimes e^{\alpha},   v\otimes e^{\beta} \rangle
=\lambda^{-\alpha^2}\delta_{\alpha,-\beta} \langle u\otimes e^0,  
v\otimes e^{0} \rangle.
\label{albet}
\end{equation}
We call this form the Li--Zamolodchikov (Li--Z) metric on ${\cal M}$.
Thus the Li--Z dual of the Fock vector $\Psi=\Psi(\mathbf{k},\mathbf{l})$ is
\[
\overline\Psi(\mathbf{k},\mathbf{l})
= \left(-1\right )^{st+\lfloor \wt(\Psi)\rfloor} 
\lambda^{2\wt(\Psi)}   
\Psi(\mathbf{l},\mathbf{k}),
\]
 where $\lfloor x\rfloor$ denotes the integer part of $x$ \cite{TZ2}. Applying \eqref{eq:lamb_eps} and \eqref{albet} it follows that 
$\Psi_{\alpha}=\Psi_{\alpha}(\mathbf{k},\mathbf{l})$ of \eqref{eq:Focka} has Li--Z dual
\begin{eqnarray}
\overline{\Psi}_{\alpha}(\mathbf{k},\mathbf{l})
&=&
 \left(-1\right )^{st+\lfloor \wt(\Psi)\rfloor}
\lambda^{2 \wt(\Psi_{\alpha})}  
\Psi_{-\alpha}(\mathbf{l},\mathbf{k})
\notag\\
&=&
 \left(-1\right )^{st+\lfloor \wt(\Psi)\rfloor}
e^{ i\pi B \wt(\Psi_{\alpha})} \,
\rho^{\wt(\Psi_{\alpha})} 
\Psi_{-\alpha}(\mathbf{l},\mathbf{k}).
\label{twisteddual}
\end{eqnarray} 

\section{Torus Intertwining $n$--Point Functions}
\label{Torus}
\subsection{Torus Intertwining $n$--Point Functions for $\Mcal$}
In this section we compute the torus intertwining $n$--point functions for $\Mcal$ and then for $V_{\Z}$. In the next section we will use these functions to define and compute the partition and $n$--point functions on a genus two Riemann surface in the $\rho$ formalism. 
Define the \lq square-bracket\rq\  vertex operator \cite{Z1}
\begin{equation*}
\Ycal[u\otimes e^{\alpha},z]=\Ycal(q_{z}^{L(0)}(u\otimes e^{\alpha}),q_{z}-1),
\label{Ysquare}
\end{equation*}%
where $q_z=e^z$. 
Define the linear operator $
T =\exp\left(\sum_{i \ge 1} r_i L(i)\right)$
associated with the exponential map 
$z\rightarrow q_z-1=\exp\left(\sum_{i\ge 1} r_i z^{i+1}\der{z}\right)z$
for $ r_1=\frac{1}{2}, r_2=-\frac{1}{12}, r_3=\frac{1}{48},\ldots$  
Then the bracket square vertex operator is given by 
\[
\Ycal[u\otimes e^{\alpha},z] = T  \Ycal\left (T ^{-1} (u\otimes e^{\alpha}),z\right)T ^{-1}. 
\]
Furthermore, as in  \cite{Z1} we find that  $(V,\Ycal[\ ,\ ],\mathbf{1},\widetilde{\omega})$ and $(V,\Ycal(\ ,\ ),\mathbf{1},\omega )$ are isomorphic generalized VOAs with $\wtil=\omega-\frac{1}{24}\vac$, a new conformal vector with vertex operator
\begin{align*}
Y[\widetilde{\omega},z] = \sum_{n\in\Z}L[n]z^{-n-2}.
\end{align*}
We let $\wt[u\otimes e^\alpha]=\wt[u]+\half\alpha^2$ denote the weight of an $L[0]$ homogeneous vector $u\otimes e^\alpha$.
\medskip

Let $M\otimes e^{\alpha}$ be an irreducible $M$-module for some $\alpha\in \mathbb{C}$ with torus partition function 
\[
Z^{(1)}_{\alpha}(q)=\Tr_{M\otimes e^{\alpha}} q^{L(0)-1/24} =\frac{q^{\half \alpha^2}}{\eta(q)},
\]
where $\eta(q)=q^{1/24}\prod_{n\ge 1}(1-q^n)$ is the Dedekind eta for modular parameter $q$. We also define the 1-point correlation function for $u\in M$ on  $M\otimes e^{\alpha}$ by
\[
Z^{(1)}_{\alpha}(u;q)=\Tr_{M\otimes e^{\alpha}}\left(o(u)q^{L(0)-1/24}\right),
\]
where $o(u)=u(1-\wt(u))$ is the Virasoro level preserving ``zero mode'' for $u$ of weight $\wt(u)$. 
In general we define the genus one intertwining $n$--point correlation function on $M\otimes e^{\alpha}$ for $n$ vectors 
$u_{1}\otimes e^{\beta_1}, \ldots , v_n\otimes e^{\beta_n}\in \Mcal$ by
\begin{eqnarray}
&&Z^{(1)}_{\alpha}  
\left( u_{1}\otimes e^{\beta_1}, z_1;  \ldots ; u_{n}\otimes e^{\beta_n}, z_n ; q\right)\notag
\\
&&=
\Tr_{M\otimes e^{\alpha}}
\left(
\Ycal \left(q_{1}^{L(0)}(u_{1}\otimes e^{\beta_1}), q_{1}\right) 
\ldots 
\Ycal \left(q_{n}^{L(0)}(u_{n}\otimes e^{\beta_n}), q_{n}\right) 
q^{L(0)-1/24}
\right),\notag  \\
\label{Znpt}
\end{eqnarray}
for formal  $q_{i}=e^{z_{i}}$ with $i=1,\ldots,n$.  
Since $e^{\beta \widehat{q}}M\otimes e^{\alpha}=M\otimes e^{\alpha+\beta}$ it follows that the $n$--point function vanishes  when $\sum_{i=1}\beta_i\neq 0$.

We next describe a natural generalization of previous results in \cite{MT1} and \cite{MTZ}. Firstly, consider the $n$--point functions for $n$ highest weight vectors $\vac\otimes e^{\beta_i}$, which we abbreviate below to $e^{\beta_i}$, for $i=1,\ldots, n$.
\begin{proposition} 
\label{prop_ZalphaK}
For $\sum_{i=1}^{n}\beta_i= 0$ then
\begin{eqnarray}
Z^{(1)}_{\alpha}  
\left( e^{\beta_1}, z_1;  \ldots ;e^{\beta_n}, z_n  ;q\right)=
\frac{q^{\half\alpha^2}}{\eta (\tau )}
\exp\left(\alpha \sum_{i=1}^{n}\beta_{i}z_{i}\right)
\prod_{1\leq r<s\leq n}K(z_{rs},\tau )^{\beta _{r}\beta _{s}},\;
\label{ZNexpalpha}
\end{eqnarray}
where $z_{rs}=z_{r}-z_{s}$ and $K(z,\tau )$ is the genus one prime form \eqref{Ktheta}.
\end{proposition}

The proof appears in Appendix~\ref{ProofZalphaK} and is a natural generalization of results  developed in \cite{MT1} and \cite{MTZ}.
\begin{remark}
\label{rem_Zdiff}
Let $F_{\alpha}^{\boldsymbol{\beta}}(\textbf{z})=Z^{(1)}_{\alpha}   
\left( e^{\beta_1}, z_1;  \ldots ;e^{\beta_n}, z_n ;q\right)$ of \eqref{ZNexpalpha}.
\begin{enumerate}
	\item[(i)]  $F_{\alpha}^{\boldsymbol{\beta}}(\textbf{z})$ is a function of the differences $z_{ij}=z_i-z_j$. Thus we may choose $z_n=0$ wlog.
 \item[(ii)] $F_{\alpha}^{\boldsymbol{\beta}}(\textbf{z})$ is analytic on $
\{z_i\in\C,\tau\in \HH: \vert z_1\vert>\ldots >\vert z_n\vert\}$. 
\item[(iii)] Note that 
\begin{eqnarray*}
&&
(z_1-z_2)^{-\beta_1\beta_2}
Z^{(1)}_{\alpha}  
\left( e^{\beta_1}, z_1;  e^{\beta_2}, z_2; \ldots ;e^{\beta_n}, z_n ;q\right)
\\
&=&
(z_2-z_1)^{-\beta_1\beta_2}
Z^{(1)}_{\alpha}  
\left( e^{\beta_2}, z_2; e^{\beta_1}, z_1;  \ldots ;e^{\beta_n}, z_n ;q\right),
\end{eqnarray*}
in compliance with \eqref{genloc} and \eqref{genlocB} (for $N=0$). Similar identities hold under all label permutations. Thus $F_{\alpha}^{\boldsymbol{\beta}}(\textbf{z})$ can be extended to all $z_i\neq z_j$ for $i\neq j$ up to branch choices.
\end{enumerate}
\end{remark}

We may similarly apply the arguments mutatis mutandis of Proposition~1 of \cite{MT1}  to obtain a closed form for the general $n$--point function \eqref{Znpt}. In particular, due to \eqref{genComm}, we may apply standard genus one Zhu recursion theory  \cite{Z1} to reduce \eqref{Znpt} to an explicit multiple of \eqref{ZNexpalpha} to find
\begin{proposition}
\label{prop_Znalpha}
For $\sum_{i=1}^{n}\beta_i= 0$ then
\begin{eqnarray*}
&&Z^{(1)}_{\alpha}  
\left( u_1\otimes e^{\beta_1}, z_1;  \ldots ;u_n\otimes e^{\beta_n}, z_n  ;q\right)\\
&&=
Q_{\alpha}^{\beta_1,\ldots,\beta_n}(u_1,z_1;\ldots u_n,z_n;q)
Z^{(1)}_{\alpha}  
\left( e^{\beta_1}, z_1;  \ldots ;e^{\beta_n}, z_n  ;q\right),
\end{eqnarray*}
where $Q_{\alpha}^{\beta_1,\ldots,\beta_n}(u_1,z_1;\ldots u_n,z_n;q)$ is an explicit  sum of elliptic and quasi-modular forms (see \cite{MT1} for details).
\end{proposition}

\subsection{Torus Intertwined $n$--Point Functions for $V_{\Z}$}
\label{genus_one}
Let $g_1,h_1$ be commuting automorphisms of  $V_{\Z}$ defined by 
\begin{equation}
\label{f_one_g}
\sigma f_1=e^{2\pi i \beta_1 a(0)}, \quad 
\sigma g_1=e^{-2\pi i \alpha_1 a(0)}, 
\end{equation} 
where $\sigma=e^{\pi i a(0)}$ is the fermion number automorphism. 
We assume that $\alpha_1,\beta_1\in \R$ so that 
\[
\theta_1=-e^{-2 \pi i \beta_1},\quad {\phi_1}= -e^{2 \pi i \alpha_1},
\]
are of unit modulus (to  ensure convergence of  Szeg\H{o} kernels below).

Following Remark~\ref{rem_Zdiff}~(i) we consider torus orbifold  intertwining $n+2$--point functions for the $\sigma g_1$--twisted module $V_{\Z+\alpha_1}=\oplus_{m\in\Z}M\otimes e^{m+\alpha_1}$ for $V_{\Z}$ defined by:
\begin{eqnarray}
&&Z^{(1)}_{V_{\Z}} \begin{bmatrix}f_1\\ g_1 \end{bmatrix} 
\left( u_{1}\otimes e^{m_1}, z_1;  \ldots , u_{n}\otimes e^{m_n}, z_n ;v_{1}\otimes e^{n_1+\kappa},w;v_{2}\otimes e^{n_2-\kappa},0; q\right)\notag
\\
\notag
&=&
\Tr_{V_{\Z+\alpha_1}}
\left(
\sigma f_1
Y \left(q_{1}^{L(0)}(u_{1}\otimes e^{m_1}), q_{1}\right) 
\ldots 
Y\left(q_{n}^{L(0)}(u_{n}\otimes e^{m_n}), q_{n}\right) 
\right).\\
&&\quad \qquad \left.
\Ycal  \left( q^{L (0)}_{w} (v_{1}\otimes e^{n_1+\kappa}), q_{w}\right) 
\Ycal  \left(v_{2}\otimes e^{n_2-\kappa}, 1\right)
q^{L(0)-1/24}
\right)
\notag
\\
&=& \sum_{\mu\in \Z+\alpha_1} e^{\tpi \mu\beta_1} Z^{(1)}_{\mu}\left(u_{1}\otimes e^{m_1}, z_1;  \ldots ; u_{n}\otimes e^{m_n}, z_n ;\right.
\notag\\
&&\quad\qquad \left.
v_{1}\otimes e^{n_1+\kappa},w;\;v_{2}\otimes e^{n_2-\kappa},0; q\right),
\label{Zkapnpt}
\end{eqnarray}
where $m_i,n_1,n_2\in\Z$ for $\kappa\in (-\half,\half)$ and  $q_{w}=e^{w}$.
The points $w,0$ will later be identified with the centres of the sewn annuli  in the genus two $\rho$-sewing scheme described in Section~\ref{Sect_Rho_g}.
 
In Propositions~14 and 15 of \cite{MTZ} all torus orbifold $n$--point functions for vectors in $V_{\Z}$ are computed by means of a generating function. 
In that analysis we made use of an explicit Zhu reduction formula developed for a VOSA with real grading. 
However, in the present case we do not have an intertwining Zhu reduction formula because of the absence of a commutator formula generalizing \eqref{Comm} for interwiners \cite{DL}. Instead, here we adopt an alternative approach for computing the intertwining $n$--point functions by exploiting the bosonized formalism. 
Thus for example,  we find 
\begin{eqnarray}
 Z^{(1)}_{V_{\Z}} \begin{bmatrix}f_1\\ g_1 \end{bmatrix} \left(e^{\kappa},w;e^{-\kappa},0; q\right) 
&=& \sum_{\mu\in \Z+\alpha_1} e^{\tpi \mu\beta_1} Z^{(1)}_{\mu}\left(e^{\kappa},w;e^{-\kappa},0; q\right)
\notag
\\
&=& \frac{1}{\eta(q)}
\frac{
\vartheta 
\begin{bmatrix}
\alpha_1 \\ \beta_1  
\end{bmatrix}
 \left(\kappa w ,  \tau \right)
}{ K (w, \tau)^{\kappa^2} },  
\label{eq:Zgh}
\end{eqnarray}
for genus one theta series \eqref{theta}.
More generally, we define a generating function for all $n+2$--point functions  \eqref{Zkapnpt} by the following formal differential form 
\begin{eqnarray}
&  
{\mathcal G}^{(1)}_n \begin{bmatrix}f_1\\ g_1 \end{bmatrix}(x_1, y_1, \ldots, x_n, y_n)&   
\notag
\\
 &   
\equiv  
Z^{(1)}_{V_{\Z}}\begin{bmatrix}f_1\\ g_1 \end{bmatrix}
\left( \psi^{+}, x_1; \; \psi^{-}, y_1; \;  \ldots ; \psi^{+}, x_n; 
\;  \psi^{-}, y_n;e^{\kappa}, w; \;  e^{-\kappa}, 0; q\right)\prod_{i=1}^n dx_i^{\half}dy_i^{\half}, &
\notag\\
\label{gen_function_0_0}
\end{eqnarray}
for  $V_{\Z}$ generators $\psi^{\pm}=e^{\pm 1}$   
alternatively inserted at  $x_i$, $y_i$ for $i=1, \ldots,  n$. 

\begin{proposition}
\label{Theorem_generating_s}
The generating form \eqref{gen_function_0_0} is given by 
\begin{eqnarray}
\label{Z2_1bos_s}
{\cal G}_n^{(1)}\begin{bmatrix}f_1\\ g_1 \end{bmatrix}(x_1, y_1, \ldots , x_n, y_n) 
= 
Z^{(1)}_{V_{\Z}}\begin{bmatrix}f_1\\ g_1 \end{bmatrix}
\left(  e^{\kappa}, w;   e^{-\kappa}, 0; q 
\right)
\; \det S_{\kappa}, 
\end{eqnarray}
where $S_{\kappa}$ denotes the $n\times n$ matrix with components 
$S_{\kappa}\left(x_i, y_j \right)$ for $i,j=1, \ldots, n$
for Szeg\H{o} kernel  \eqref{skappa_1}.  
\end{proposition}
\noindent{\textbf{Proof.}}
 Proposition~\ref{prop_ZalphaK} implies that 
\begin{eqnarray*}
&&{\mathcal G}^{(1)}_n \begin{bmatrix}f_1\\ g_1 \end{bmatrix}(x_1, y_1, \ldots, x_n, y_n)
\\
&=& \frac{1}{\eta(\tau)K (w, \tau)^{\kappa^2} } 
 \vartheta  \begin{bmatrix}\alpha_1\\ \beta_1 \end{bmatrix}
\left(\sum\limits_{m=1}^n (x_m - y_m) +\kappa w, \tau\right)
\\
& & 
  \cdot 
\frac{ 
 \prod\limits_{1 \le i<j \le n}  K (x_i - x_j, \tau) K (y_i - y_j, \tau) 
}
{ 
\prod\limits_{1 \le i, j \le n} K (x_i - y_j, \tau) 
}
\prod\limits_{1 \le i, j \le n}
\left[\frac{  
K (x_i- w, \tau) K (y_j, \tau)
}
{ 
 K (x_i, \tau)  K (y_j -w, \tau) 
}
\right]^\kappa dx_i^{\half}dy_j^{\half}.   
\end{eqnarray*}
Next we recall the Frobenius identity (e.g. \cite{F}, \cite{MTZ})
\begin{equation*}
\label{Fay_trisecant}
\frac{\vartheta  \begin{bmatrix}\alpha_1\\ \beta_1 \end{bmatrix} 
\left( \sum\limits_{m=1}^{n}(x_{m}-y_{m}), \tau \right)  }
{
\vartheta  \begin{bmatrix}\alpha_1\\ \beta_1 \end{bmatrix} 
\left( 0, \tau \right) }
\frac{ \prod\limits_{1\leq i<j\leq n} K (x_{i}-x_{j},\tau) \; K (y_{i}-y_{j},\tau )}
{\prod\limits_{1\leq i,j\leq n}\; K (x_{i}-y_{j},\tau ) }=
\det P_1, 
\end{equation*}
where $P_1$ denotes the $n\times n$ matrix with twisted Weierstrass function components 
$P_1
\left[\begin{smallmatrix}\theta_1\\ \phi_1 \end{smallmatrix}\right] 
(x_i,y_j)$ of    \eqref{p_one_theta} for $i,j=1, \ldots, n$.  
Noting that 
$
\vartheta 
\left[\begin{smallmatrix}\alpha_1\\ \beta_1 \end{smallmatrix}\right] 
\left( x+\kappa w, \tau \right)
=\vartheta  \left[\begin{smallmatrix}\alpha_1\\ \beta_1-\frac{\kappa w}{\tpi}  \end{smallmatrix}\right] 
\left( x, \tau \right)$
we obtain
\begin{eqnarray*}
&&
 \vartheta \begin{bmatrix}\alpha_1\\ \beta_1 \end{bmatrix}
\left( \sum\limits_{m=1}^{n}(x_{m}-y_{m})+\kappa w, \tau \right) 
\frac{ \prod\limits_{1\leq i<j\leq n} K (x_{i}-x_{j},\tau) \; K (y_{i}-y_{j},\tau )}
{\prod\limits_{1\leq i,j\leq n}\; K (x_{i}-y_{j},\tau ) }\\
&=&
\vartheta \begin{bmatrix}\alpha_1\\ \beta_1 \end{bmatrix}
\left( \kappa w, \tau \right)
\det \left[
\frac{\vartheta 
\begin{bmatrix}
\alpha_1 \\ 
\beta_1
\end{bmatrix} (x_{i}-y_{j}+\kappa w,\tau )}
{\vartheta 
\begin{bmatrix}
\alpha_1 
\\ 
\beta_1
\end{bmatrix}(\kappa w,\tau )}\frac{1}{K (x_{i}-y_{j},\tau)}
\right].
\end{eqnarray*}
The result follows on absorbing the $\left[\frac{  
K (x_i- w, \tau) K (y_j, \tau)
}
{ 
 K (x_i, \tau)  K (y_j -w, \tau) 
}
\right]^\kappa dx_i^{\half}dy_j^{\half}$ factors into the determinant and using \eqref{eq:Zgh}. $\square$

\medskip

Finally, we obtain the following generalization of Proposition~15 of \cite{MTZ} concerning the generating properties of \eqref{gen_function_0_0}.
\begin{proposition}
\label{prop:Ggen}
$ {\mathcal G}^{(1)}_m
 \left[\begin{smallmatrix}f_1\\ g_1 \end{smallmatrix}\right]
(x_1, y_1, \ldots, x_m, y_m)$ is a generating function for all 
torus orbifold intertwining $n+2$--point functions \eqref{Zkapnpt}. 
In particular, for a pair of square bracket mode twisted Fock vectors \eqref{eq:Focka} 
\begin{eqnarray*}
\Psi_{\kappa}[\textbf{k}_1,\textbf{l}_2]
&=&
e^{\kappa \qhat}\psi^{+}[-k_{11}]\ldots\psi^{+}[-k_{1s_1}]
\psi^{-}[-l_{21}]\ldots\psi^{-}[-l_{2t_2}]\vac\\
\Psi_{-\kappa}[\textbf{k}_2,\textbf{l}_1]
&=&
e^{-\kappa \qhat}\psi^{+}[-k_{21}]\ldots\psi^{+}[-k_{2s_2}]
\psi^{-}[-l_{11}]\ldots\psi^{-}[-l_{1t_1}]\vac,
\end{eqnarray*}
for $p=s_1+s_2=t_1+t_2>0$ then
\begin{eqnarray}
\label{ZPsi2}
&&Z^{(1)}_{V_{\Z}} \begin{bmatrix}f_1\\ g_1 \end{bmatrix} 
\left(\Psi_{\kappa}[\textbf{k}_1,\textbf{l}_2],w;
\Psi_{-\kappa}[\textbf{k}_2,\textbf{l}_1],0; q\right)
\notag\\ 
&&=\epsilon\;
Z^{(1)}_{V_{\Z}} \begin{bmatrix}f_1\\ g_1 \end{bmatrix}\left( 
 e^{\kappa}, w; \; 
 e^{-\kappa}, 0; q \right)
 \det C_{ab}(\textbf{k}_a,\textbf{l}_b)
\end{eqnarray}
where 
\[
\epsilon=(-1)^{(t_1+s_2)t_2+\lfloor \half p\rfloor}e^{i\pi B\kappa(s_2-t_1)},
\] 
 for odd integer $B$ of 
\eqref{genlocB} and 
\[
C_{ab}(\textbf{k}_a,\textbf{l}_b)
=\begin{bmatrix}
C_{11}(\textbf{k}_1,\textbf{l}_1)& C_{12}(\textbf{k}_1,\textbf{l}_2)\\
C_{21}(\textbf{k}_2,\textbf{l}_1) & C_{22}(\textbf{k}_2,\textbf{l}_2)
\end{bmatrix}
\]
 is a $p\times p$ block matrix with components 
$
C_{ab}(k_{ai_a},l_{bj_b})
$
for $i_a=1,\ldots,s_a$ and  $j_b=1,\ldots ,t_b$ for $a,b=1,2$ for $S_{\kappa}$ expansion
coefficients  \eqref{Cijdef}.
\end{proposition}

\noindent{\textbf{Proof.}}
Consider the generating function for $\psi^{+}$ and $\psi^{-}$  alternatively inserted  at
$
x_{21}+w,\ldots,x_{2s_{2}}+w,x_{11},\ldots,x_{1s_{1}}$ and 
$y_{21}+w,\ldots,y_{2t_{2}}+w,y_{11},\ldots,y_{2t_{1}}$
respectively.
Reordering the vertex operators and using generalized locality 
\eqref{genlocB} we obtain
\begin{eqnarray}
&&
\widehat{\epsilon}\,
Z^{(1)}_{V_{\Z}} \begin{bmatrix}f_1\\ g_1 \end{bmatrix} 
\left( \psi^{+}, x_{21}+w; \ldots; \psi^{+}, x_{2s_{2}}+w; 
\psi^{-}, y_{21}+w;\ldots ;\psi^{-}, y_{2t_{2}}+w;
e^{\kappa},w;\right.
\notag
\\
&&
\qquad\qquad
\left.
 \psi^{+}, x_{11}; \ldots ;\psi^{+}, x_{1s_{1}}; 
\psi^{-}, y_{11};\ldots ;\psi^{-}, y_{1t_{1}};
e^{-\kappa},0; 
q\right)
\prod_{i_a,j_b} dx_{ai_{a}}^{\half}dy_{bj_{b}}^{\half},\notag
\\
\label{epsZ}
\end{eqnarray}
for $\widehat{\epsilon}=(-1)^{s_2 t_2+\lfloor \half p\rfloor}e^{-i\pi B\kappa(s_2-t_1)}$.
Generalized associativity \eqref{genassoc} implies that 
\begin{eqnarray*}
&&
\prod_{i,j}^{s_1}
(x_{2i}+w)^{M_i-\kappa}
(y_{2j}+w)^{N_j+\kappa}
\Ycal\left[\psi^{+}, x_{21}+w \right] \ldots\Ycal\left[\psi^{-}, y_{2t_{2}}+w\right]
\Ycal\left[e^{\kappa},w\right]
\notag
\\
&=&
\prod_{i,j}^{s_1}
(w+x_{2i})^{M_i-\kappa}
(w+y_{2j})^{N_j+\kappa}
\Ycal\left[
\Ycal\left[\psi^{+}, x_{21} \right] \ldots\Ycal\left[\psi^{-}, y_{2t_{2}}\right]
e^{\kappa},w\right],
\end{eqnarray*}
for $M_i,N_j\gg 0$. But using \eqref{eq:qhatconj} and 
$\Delta(\kappa,z)\psi^{\pm}=z^{\pm\kappa}\psi^{\pm}$ 
we find that
\begin{eqnarray*}
&&Y[\psi^{+},x_{21}]\ldots Y[\psi^{+},x_{2s_1}] Y[\psi^{-},y_{21}]\ldots Y[\psi^{-},y_{2t_2}]e^{\kappa}\\
&=&
\prod_{i=1}^{s_1}x_{2i}^{\kappa}
\prod_{j=1}^{t_2}y_{2j}^{-\kappa}
\cdot
e^{\kappa \widehat{[q]}}
Y[\psi^{+},x_{21}]\ldots Y[\psi^{+},x_{2s_1}] Y[\psi^{-},y_{21}]\ldots Y[\psi^{-},y_{2t_2}]\vac\\
&=&
\sum_{\textbf{k}_1\in \Z^{s_1},\textbf{l}_2\in \Z^{t_2}}
\prod_{i=1}^{s_1}x_{2i}^{k_{1i}+\kappa-1}
\prod_{j=1}^{t_2}y_{2j}^{l_{2j}-\kappa-1}
\,
\Psi_{\kappa}[\textbf{k}_1,\textbf{l}_2], 
\end{eqnarray*}
generates the Fock vector $\Psi_{\kappa}[\textbf{k}_1,\textbf{l}_2]$ as the coefficient of 
$\prod_{i=1}^{s_1}x_{2i}^{k_{1i}+\kappa-1}
\prod_{j=1}^{t_2}y_{2j}^{l_{2j}-\kappa-1}$ (and where $\left[a[n],\widehat{[q]}\right]=\delta_{n0}$ is the square bracket version of \eqref{pcom1}). 
Similarly $\Psi_{-\kappa}[\textbf{k}_2,\textbf{l}_1]$ is the coefficient of
$\prod_{i=1}^{s_2}x_{1i}^{k_{2i}-\kappa-1}
\prod_{j=1}^{t_1}y_{1j}^{l_{1j}+\kappa-1}$ in
\[
Y[\psi^{+},x_{11}]\ldots Y[\psi^{+},x_{1s_2}] Y[\psi^{-},y_{11}]\ldots Y[\psi^{-},y_{1t_1}]e^{-\kappa}.
\]
\eqref{ZPsi2} follows by reading off the coefficient in the corresponding expansion  \eqref{Cijdef} of the Szeg\H{o} kernels in the determinant of  \eqref{Z2_1bos_s}. Finally, interchanging the first $t_1$ and second $t_2$ columns of the resulting determinant gives $(-1)^{t_1 t_2} \det C_{ab}(\textbf{k}_a,\textbf{l}_b)$. Thus 
$ \epsilon=(-1)^{t_1 t_2}/\widehat{\epsilon}$ as given. In general, we can  expand  \eqref{Z2_1bos_s} in appropriate parameters much as in Proposition~15 of \cite{MTZ} to obtain any $n+2$--point function for $V_{\Z}$-Fock vectors inserted at $z_1,\ldots ,z_n$ for $z_i\neq 0,w$ together with the above pair of twisted Fock vectors.
$\square$

\section{Genus Two Twisted Partition and $n$--Point Functions}
\label{genus_two_n-point}
We now arrive at the main results of this paper, namely, the  computation of the genus two twisted partition and $n$--point function for $V_{\Z}$ on the genus two Riemann surface formed by self-sewing a torus as described in Section~\ref{Sect_Rho_g}. 
Our definitions are a natural extension of those for a VOA introduced in \cite{MT3}. In the present case, we  consider genus two twisted partition and $n$--point function which depend on 
four $V_{\Z}$-automorphisms:
\begin{eqnarray}
\label{f_g}
\sigma f_1 &=& e^{2\pi i \beta_1 a(0)}, 
\quad 
\sigma g_1 = e^{-2\pi i \alpha_1 a(0)}, 
\notag
\\
\sigma f_2 &=& e^{2\pi i \beta_2 a(0)}, 
\quad 
\sigma g_2=e^{-2\pi i \kappa a(0)},\label{eq:fgi}
\end{eqnarray} 
for real $\alpha_1,\beta_1,\beta_2$ and $\kappa\in (-\half,\half)$.  
We let $f=(f_1,f_2 )$ and $g=(g_1, g_2)$ denote pairs of automorphisms. 
We also define unit modulus parameters
\begin{eqnarray}
\theta_1&=&-e^{-2 \pi i \beta_1},
\quad 
\phi_1= -e^{2 \pi i \alpha_1},
\notag
\\
\theta_2 &=&-e^{-2 \pi i \beta_2},
\quad 
\phi_2 = -e^{2 \pi i \kappa},
\label{eq:thetaphi}
\end{eqnarray}
later identified with the multipliers \eqref{eq:periodicities} of a genus two Szeg\H{o} kernel below.
    

Following \cite{MT3} we define the continuous orbifold genus two $n$--point function for 
$u_{i}\otimes e^{m_i}\in V_{\Z}$ inserted at $z_i$ on the sewn genus two Riemann surface in
the $\rho$ sewing scheme as the following infinite sum of genus one intertwined $n+2$ point functions:
\begin{eqnarray}
\notag
&& Z^{(2)}_{V_{\Z}}
\begin{bmatrix}f\\g\end{bmatrix} (u_{1}\otimes e^{m_1}, z_1;  \ldots , u_{n}\otimes e^{m_n}, z_n; \tau, w, \rho )  
 \\
&&
= 
\sum_{\Psi_{\kappa}}
Z^{(1)}_{V_{\Z}} \begin{bmatrix}f_1\\ g_1 \end{bmatrix} 
\left( u_{1}\otimes e^{m_1}, z_1;  \ldots , u_{n}\otimes e^{m_n}, z_n ;\sigma f_2\Psi_{\kappa},w;\overline{\Psi}_{\kappa},0; q\right),\quad
\label{Ztwonpt}
\end{eqnarray}
where the sum is taken over the square bracket Fock basis $\Psi_{\kappa}\in V_{\Z+\kappa}$  with square bracket Li-Z dual $\overline{\Psi}_{\kappa}$ for $\lambda=(e^{i\pi B } \rho)^{\half }$ of \eqref{eq:lamb_eps} and with $\sigma f_2$ acting on $\Psi_\kappa$.

\subsection{Genus Two Twisted Partition Function}
The partition or $0$--point function is expressible in terms of the basic genus one twisted 2-point function \eqref{eq:Zgh} and the infinite Szeg\H{o} moment matrix  
$T=T
(\tau,w,\rho\vert \theta_1,\theta_2,\phi_1,\kappa)$ of \eqref{T_matrix} as follows:
\begin{theorem}
\label{Theorem_Z2_ferm}
The  genus two continuous orbifold partition function for $V_{\Z}$  on a Riemann surface in the $\rho$-sewing formalism is convergent on ${\mathcal D}^{\rho}$ and
 is given by
\begin{eqnarray}
Z^{(2)}_{V_{\Z}}
\begin{bmatrix}f\\ g \end{bmatrix}(\tau, w, \rho)
&=&
e^{2 i\pi \beta_2\kappa }(e^{i\pi B }\rho)^{\half \kappa^2}\,
Z^{(1)}_{V_{\Z}} \begin{bmatrix}f_1\\ g_1 \end{bmatrix}\left( 
 e^{\kappa}, w; \; 
e^{-\kappa}, 0; \; q\right)
\notag
\\
&&\cdot \det \left(I - T \right). 
\label{Z2_ferm}
\end{eqnarray}
\end{theorem}

In order to prove Theorem \ref{Theorem_Z2_ferm} we first recall Proposition~2 of \cite{TZ2}. 
Let $R$ be an $P\times P$ matrix and let $\mathbf{m}=(m_{1}, \ldots, m_{p})$ 
denote $p$ ordered subindices with $1\le m_{1}<\ldots <m_{p}\le P$. We refer to $\mathbf{m}$ as 
an \emph{$P$-subindex of length $p$}. Let 
\begin{equation*}
R(\mathbf{m},\mathbf{m})=\left( R(m_{r},m_{s})\right), 
\end{equation*}
for $ r,s=1,\ldots, p$ denote the $p\times p$ submatrix of $R$. We define $R(\mathbf{m},\mathbf{m})=1$ in the degenerate case $p=0$ and let $I$ denote the identity matrix. We have
\begin{proposition}
\label{propdet} 
\begin{eqnarray}
\label{detR}
\det \left(I+R\right) = 
\sum\limits_{p =0}^{P}
\sum\limits_{\mathbf{m}}
 \det R(\mathbf{m},\mathbf{m}), 
\end{eqnarray}
where the inner sum runs over all $P$-subindices of length $p$.  
 $\square$
\end{proposition}
\noindent{\bf Proof of Theorem \ref{Theorem_Z2_ferm}.} 
Recall the Fock vector $\Psi_{\kappa}=e^{\kappa \qhat}\Psi$ 
with
\[
\Psi=
\Psi[\mathbf{k},\mathbf{l}]=\psi^{+}[-k_{1}]\ldots\psi^{+}[-k_{s}]
\psi^{-}[-l_{1}]\ldots\psi^{-}[-l_{t}]\vac,
\]
 has Li--Z dual $\overline{\Psi}_{\kappa}[\mathbf{k},\mathbf{l}]
=\epsilon_{1} \rho^{\wt[\Psi_{\kappa}]}
\Psi_{-\kappa}[\mathbf{l},\mathbf{k}]$
where 
\begin{eqnarray}
\epsilon_{1}
=
 \left(-1\right )^{st+\lfloor \wt[\Psi]\rfloor}
e^{ i\pi B \wt[\Psi_{\kappa}]},
\label{eq:eps1}
\end{eqnarray} 
for $\wt[\Psi_{\kappa}]=\wt[\Psi]+\kappa(s-t)+\half \kappa^2$. 
Furthermore, 
\[
\sigma f_2\Psi_{\kappa}=(-\theta_2)^{t-s}e^{2 i\pi \beta_2\kappa }\Psi_{\kappa},
\]
from \eqref{eq:fgi} and \eqref{eq:thetaphi}.
Thus we apply \eqref{ZPsi2} of Proposition~\eqref{prop:Ggen} with $\textbf{k}=\textbf{k}_1=\textbf{l}_1$ for $s=s_1=t_1$ and  $\textbf{l}=\textbf{k}_2=\textbf{l}_2$ for $t=s_2=t_1$ and $p=s+t$ to obtain 
\begin{eqnarray*}
\notag
Z^{(2)}_{V_{\Z}} \begin{bmatrix}f\\g \end{bmatrix} ( \tau, w, \rho )  
&=&
e^{2 i\pi \beta_2\kappa }
Z^{(1)}_{V_{\Z}} \begin{bmatrix}f_1\\ g_1 \end{bmatrix}\left( 
 e^{\kappa}, w; \; 
 e^{-\kappa}, 0; q \right)
\\
&&\cdot
\sum_{\textbf{k},\textbf{l}}
\epsilon_1 \epsilon_2(-\theta_2)^{t-s}\rho^{\wt[\Psi_{\kappa}]}
 \det 
\begin{bmatrix}
C_{11}(\textbf{k},\textbf{k})& C_{12}(\textbf{k},\textbf{l})\\
C_{21}(\textbf{l},\textbf{k}) & C_{22}(\textbf{l},\textbf{l})
\end{bmatrix},
\end{eqnarray*}
where 
\begin{equation}
\epsilon_2=(-1)^{t+st+\lfloor \half p\rfloor}e^{i\pi B\kappa(t-s)}.
\label{eq:eps2}
\end{equation}
Noting that $
 \left(-1\right )^{\lfloor \wt[\Psi]\rfloor}
e^{ i\pi B \wt[\Psi]}=(-1)^{\lfloor \half p\rfloor}e^{\half i\pi B p}$
we  find
\[
\epsilon_1 \epsilon_2(-\theta_2)^{t-s}\rho^{\wt[\Psi_{\kappa}]}
=
(-\xi\theta_2^{-1})^{s}
(\xi\theta_2)^{t}
\rho^{\wt[\Psi]+\kappa(s-t)}
(e^{i\pi B}\rho)^{\half \kappa^2},
\]
for  $\xi=e^{\half i\pi B }$ of  \eqref{dz1dz2rho}. 
We may absorb all but the last of these factors into the $p\times p$ determinant to find
\begin{eqnarray*}
&&(-\xi\theta_2^{-1})^{s}
(\xi\theta_2)^{t}
\rho^{\wt[\Psi]+\kappa(s-t)}
\det   
\begin{bmatrix}
C_{11}(k_i,k_j) & C_{12}(k_i,l_{j'})\\
C_{21}(l_{i'},k_j) & C_{22}(l_{i'},l_{j'})
\end{bmatrix}
\\
&&=\det   
\begin{bmatrix}
-\xi\theta_2^{-1}\rho^{\half(k_i+k_j+2\kappa-1)} C_{11}(k_i,k_j) 
& \xi\theta_2\,\rho^{\half (k_i+l_{j'}-1)} C_{12}(k_i,l_{j'})\\
-\xi\theta_2^{-1}\rho^{\half(l_{i'}+k_j-1)} C_{21}(l_{i'},k_j) 
& \xi\theta_2\,\rho^{\half(l_{i'}+l_{j'}-2\kappa-1)} C_{22}(l_{i'},l_{j'})
\end{bmatrix}\\
&&=\det\left(- T(\mathbf{m},\mathbf{m})\right),
\end{eqnarray*}
for $T$ of  \eqref{T_matrix} with $i,j=1,\ldots,s$ and $i',j'=1,\ldots,t$ and where
$\mathbf{m}=(\mathbf{k},\mathbf{l})$ i.e. 
$m_{i}=k_{i}$ and $m_{i'+s}=l_{i'}$. But $\mathbf{m}$ is a
$P$-subindex of length $p$ for some sufficiently large $P$ so that we may apply Proposition~\ref{propdet} for sufficiently large $P$ by truncating the partition function to any finite order in $\rho$ to obtain the result. 
Finally, convergence on ${\mathcal D}^\rho$ is guaranteed by Theorem~\ref{theorem_S2}. $\square$

\subsection{The Genus Two $n$--Point Generating Form}
\label{gen_n_point_functions}
Similarly to the genus one situation 
\eqref{gen_function_0_0},
we define a genus two continuous orbifold generating differential form by 
\begin{eqnarray}
\notag
& & 
{\mathcal G}^{(2)}_n\begin{bmatrix}f\\g \end{bmatrix}
(x_1, y_1,  \ldots, x_n, y_{n})
\\
 & &  
 = 
Z^{(2)}_{V_{\Z}} \begin{bmatrix}f\\g \end{bmatrix} (\psi^{+}, x_1; \psi^{-}, y_1; \ldots , 
\psi^{+}, x_n; \psi^{-}, y_n; \tau, w, \rho )
\prod_{i=1}^n dx_i^{\half}dy_i^{\half}.
\notag\\
\label{gen_function_2}  
\end{eqnarray}
Recalling  \eqref{skappa_1}, \eqref{hdef1}, \eqref{hdef2} and \eqref{s_genus_two} we define matrices
\begin{eqnarray}
\label{stwo}
& S^{(2)}=\left(S^{(2)}(x_{i}, y_{j})\right),
\quad
S_{\kappa} =\left(S_{\kappa}(x_{i},y_{j})\right),&
\\
\notag
& H=\left( \left(h(x_{i})\right) (k,a) \right), \quad
\overline{H}^{T}=\left( \left( \overline{h} (y_{i}) \right) (l,b) \right)^{T}.& 
\end{eqnarray}
$S^{(2)}$ and  $S_{\kappa}$ are finite matrices indexed by 
$i,j=1,\ldots, n$; $H$  is semi-infinite with $n$ rows indexed by $i$ and columns
 indexed by $k\ge 1$ and $a=1,2$ and $\overline{H}^{T}$  
is semi-infinite with rows indexed by $l\ge 1$ and $b=1,2$ and with $n$ columns indexed by $j$.

We can now state the second main result of this paper.  
\begin{theorem} 
\label{generating_n_point_rank_two}
The generating form \eqref{gen_function_2} is given by 
\begin{eqnarray}
& & {\mathcal G}^{(2)}_n \begin{bmatrix}f\\g \end{bmatrix}
(x_1, y_1,  \ldots, x_n, y_{n})
=
 Z^{(2)}_{V_{\Z}}
\begin{bmatrix}f\\g \end{bmatrix}
\left(\tau, w, \rho \right) \;
 \det 
 S^{(2)}, \notag
\\
\label{Z2_def_epsorbi_two_point_1_alpha}
\end{eqnarray}
for $n\times n$ matrix  $S^{(2)}$ of \eqref{stwo} and  $Z^{(2)}_{V_{\Z}}
\left[\begin{smallmatrix}f\\g \end{smallmatrix}\right] 
\left(\tau, w, \rho \right)$ is the genus two twisted partition function \eqref{Z2_ferm}.  
\end{theorem}
Normalized relative to the genus two partition function, the  2-point function for $\psi^{+}$ and $\psi^{-}$ is thus  given by the Szeg\H{o} kernel and more generally, the generating function is a Szeg\H{o} kernel determinant. This agrees with the \emph{assumed} form in \cite{R} or as found by string theory methods using a Schottky parametrization in \cite{DVPFHLS} and in the genus two $\epsilon$-sewing
scheme \cite{TZ2}. 

\medskip

In order to prove Theorem \ref{generating_n_point_rank_two} we firstly note that similarly to  Proposition~3 in  \cite{TZ2} we have 
\begin{proposition}
\label{prop_Szegodet}
\begin{equation}
\label{det_product}
\det \left[ 
\begin{array}{cc}
S_{\kappa} & -\xi  HD^{\theta_2}   
\\
\overline{H}^{T} & I - T
\end{array}
\right] 
=\det  S^{(2)} \; \det(I - T), 
\end{equation}
with  $T$, $D^{\theta_2}$ of \eqref{T_matrix} and \eqref{Dtheta}.  
 $\square$
 \end{proposition}
We also note the following generalization of Proposition~\ref{propdet}. 
Let $S,U,V,R$ be $n\times n$, $n\times P$, $P\times n$ and $P\times P$ matrices respectively. 
For $\mathbf{m}$, a $P$-subindex of length $p$, we let 
\begin{eqnarray*}
U(\mathbf{m})&=&\left( U(i,m_{s})\right),\quad i=1,\ldots, n;\; s=1,\ldots p, \\
V(\mathbf{m})&=&\left( V(m_{r},j)\right),\quad \; r=1,\ldots p;\; j=1,\ldots, n.
\end{eqnarray*}
We then find 
\begin{proposition}
\label{propdet2} 
\begin{eqnarray}
\label{detR2}
\det 
\begin{bmatrix}
S & U\\
V & I+R
\end{bmatrix}
 = 
\sum\limits_{p =0}^{P}
\sum\limits_{\mathbf{m}}
 \det 
\begin{bmatrix}
S & U(\mathbf{m})\\
V(\mathbf{m}) & R(\mathbf{m},\mathbf{m})
\end{bmatrix}, 
\end{eqnarray}
where the inner sum runs over all $P$-subindices of length $p$.  
 $\square$
\end{proposition}
\noindent
{\bf Proof of Theorem \ref{generating_n_point_rank_two}.}  
Repeating the notation of the proof of Theorem~\ref{Theorem_Z2_ferm} we find that
$G^{(2)}_n(x_1,\ldots,y_n)\prod_{i=1}^n dx_i^{-\half}dy_i^{-\half}$ is 
\begin{eqnarray*}
e^{2 i\pi \beta_2\kappa }
\sum_{\textbf{k},\textbf{l}}
\epsilon_1 (-\theta_2)^{t-s}
\rho^{\wt[\Psi_{\kappa}]}
Z^{(1)}_{V_{\Z}} 
\begin{bmatrix}f_1\\ g_1 \end{bmatrix}
\left(\psi^{+},x_1;\ldots;\psi^{-},y_{n}; 
 \Psi_{\kappa}[\textbf{k},\textbf{l}], w; 
\Psi_{-\kappa}[\textbf{l},\textbf{k}], 0; q \right),
\end{eqnarray*}
with $\epsilon_1$ of \eqref{eq:eps1}.
The expansion technique described in Proposition 
\ref{prop:Ggen} implies  
\begin{eqnarray*}
&&Z^{(1)}_{V_{\Z}} 
\begin{bmatrix}f_1\\ g_1 \end{bmatrix}\left(\psi^{+},x_1;\ldots;\psi^{-},y_{n}; 
 \Psi_{\kappa}[\textbf{k},\textbf{l}], w;  
\Psi_{-\kappa}[\textbf{l},\textbf{k}], 0; q \right)\prod_{i=1}^n dx_i^{\half}dy_i^{\half}\\
&=&
\epsilon_2\;
Z^{(1)}_{V_{\Z}} 
\begin{bmatrix}f_1\\ g_1 \end{bmatrix}\left( 
 e^{\kappa}, w; 
 e^{-\kappa}, 0; q \right)
\det 
\begin{bmatrix}
S_{\kappa}(x_{i},y_{j}) & d_{1}(x_{i},\textbf{k}) & d_{2}(x_{i},\textbf{l})\\
\overline{d}_{1}(y_{j},\textbf{k}) & C_{11}(\textbf{k},\textbf{k})& C_{12}(\textbf{k},\textbf{l})\\
\overline{d}_{2}(y_{j},\textbf{l}) & C_{21}(\textbf{l},\textbf{k}) & C_{22}(\textbf{l},\textbf{l})
\end{bmatrix},
\end{eqnarray*}
for  $\epsilon_2$ of \eqref{eq:eps2} and $d_{a}(x,k),\overline{d}_a(y,k)$ of \eqref{ddef1} and \eqref{ddef2}. 
Absorbing various factors as before we find
\begin{eqnarray*}
&&(-\xi )^p (-1)^{t}\theta_2^{t-s}\rho^{\wt[\Psi]+\kappa(s-t)}
\det 
\begin{bmatrix}
S_{\kappa}(x_{i},y_{j}) & d_{1}(x_{i},\textbf{k}) & d_{2}(x_{i},\textbf{l})\\
\overline{d}_{1}(y_{j},\textbf{k}) & C_{11}(\textbf{k},\textbf{k})& C_{12}(\textbf{k},\textbf{l})\\
\overline{d}_{2}(y_{j},\textbf{l}) & C_{21}(\textbf{l},\textbf{k}) & C_{22}(\textbf{l},\textbf{l})
\end{bmatrix}
\\
&&= \det 
\begin{bmatrix}
S_{\kappa} & -\xi \left( HD^{\theta_2}\right)(\mathbf{m})\\
\left(\overline{H}^{T}\right)(\mathbf{m}) & -T(\mathbf{m},\mathbf{m})
\end{bmatrix},
\end{eqnarray*}
for $\mathbf{m}=(\mathbf{k},\mathbf{l})$ as before.
The result follows on applying Proposition~\ref{propdet2} and followed by Proposition~\ref{prop_Szegodet}. $\square$

\subsection{Bosonization}
\label{boson_prod_form}
We may compute the genus two twisted partition and generating functions in an alternative way by use of  a bosonic basis $u\otimes e^{\nu}\in V_{\Z+\kappa}$ for $u\in M$ and $\nu \in\Z+\kappa$.   
In particular, we can immediately exploit and extend the results for lattice VOAs in  \cite{MT3} in the light of Section~\ref{Torus} above. 

Let us firstly recall from \cite{MT2, MT3} the following definitions: 
\begin{eqnarray}
P_{2}(\tau ,z) &=&\wp (\tau ,z)+E_{2}(\tau )  \notag \\
&=&\frac{1}{z^{2}}+\sum_{k=2}^{\infty }(k-1)E_{k}(\tau )z^{k-2},  \label{P2}
\end{eqnarray}%
for  Weierstrass function $\wp (\tau ,z)$ and Eisenstein series $E_{k}(\tau )$. We define
$P_{k+1}(\tau ,z)=-\frac{1}{k}\partial_z P_{k}(\tau ,z)$ for $k\ge 2$  and for 
 $k,l\geq 1$ 
\begin{eqnarray*}
C(k,l,\tau) &=&(-1)^{k+1}\frac{(k+l-1)!}{(k-1)!(l-1)!}E_{k+l}(\tau ),
\label{Ckldef} \\
D(k,l,\tau ,z) &=&(-1)^{k+1}\frac{(k+l-1)!}{(k-1)!(l-1)!}%
P_{k+l}(\tau ,z),  \label{Dkldef}
\\
R_{ab}(k,l)&=&-\frac{\rho ^{(k+l)/2}}{\sqrt{kl}}
\begin{bmatrix}
D(k,l,\tau ,w) & C(k,l,\tau ) \\ 
C(k,l,\tau ) & D(l,k,\tau ,w)%
\end{bmatrix}.\label{Rkldef}
\end{eqnarray*}
Theorems~5.1 and 5.6 of \cite{MT3} tell us that:
\begin{theorem}
\label{Theorem_Z2_boson_rho}
The genus two partition function for the Heisenberg VOA $M$ in the $\rho$ sewing scheme is 
\begin{equation}
Z_{M}^{(2)}(\tau ,w,\rho )=Z_{M}^{(1)}(\tau ) \det (1-R)^{-1/2},
\label{Z2_1bos_rho}
\end{equation}
where $Z_{M}^{(1)}(\tau )=1/\eta(q)$ is the genus one partition function. $Z_{M}^{(2)}(\tau ,w,\rho )$ is holomorphic on ${\mathcal D}^{\rho}$.
\end{theorem}
Theorem~6.1 of \cite{MT3} concerns the genus two partition function for a lattice VOA. We obtain a natural generalization for a twisted genus two partition function defined by
\begin{eqnarray}
&& Z^{(2)}_{\mu,\nu} ( \tau, w, \rho )  
= 
\sum_{\Psi_{\nu}\in M\otimes e^{\nu}}
Z^{(1)}_{\mu} 
\left( \Psi_{\nu},w;\overline{\Psi}_{\nu},0; q\right),
\label{Ztwomunu}
\end{eqnarray}
for $\mu,\nu\in \C$ where the sum is taken over a $ M\otimes e^{\nu}$ basis, $\overline{\Psi}_{\nu}$ is the square bracket  Li-Z dual and the summand is a genus one 2-point function \eqref{Znpt} for $\Mcal$. We  find that 
\begin{theorem}
\label{Theorem_Z2_munu}
\begin{equation}
Z^{(2)}_{\mu,\nu} ( \tau, w, \rho ) =
e^{i\pi(\mu^2\Omega_{11}+2\mu\nu\Omega_{12}+\nu^2\Omega_{22})}
Z_{M}^{(2)}(\tau,w,\rho ),
\label{Z2_1bos_munurho}
\end{equation}
where $\Omega$ is the genus two period matrix.
\end{theorem}
\textbf{Proof.}  The proof follows precisely that for Theorem~6.1 of \cite{MT3} by use of a bosonic Fock basis and by applying
Propositions~\ref{prop_ZalphaK} and \ref{prop_Znalpha}. $\square$ 
\medskip

We are now able to compute the genus two twisted partition function \eqref{Z2_ferm} by use of a bosonic basis to obtain
\begin{theorem}
\label{Theorem_Z2_fermboson}
The genus two twisted partition function for $V_{\Z}$ on a Riemann surface in the $\rho$-sewing formalism 
 is given by
\begin{eqnarray}
Z^{(2)}_{V_{\Z}}
\begin{bmatrix}f\\g \end{bmatrix}(\tau, w, \rho)
=\vartheta
 \begin{bmatrix}
\alpha \\ \beta  
\end{bmatrix}
(  \Omega)\;
Z^{(2)}_{M}(\tau, w, \rho), \quad
\label{Z2_fermboson}
\end{eqnarray}
for genus two Riemann theta function \eqref{theta} with characteristics $\alpha_1,\alpha_2=\kappa,\beta_1,\beta_2$ of \eqref{eq:thetaphi} and  where $Z^{(2)}_{M}(\tau, w, \rho)$ is the genus two Heisenberg partition function. 
\end{theorem}
\textbf{Proof.}
Let $\Psi_{\nu}=u\otimes e^{\nu}\in V_{\Z+\kappa}$ with  square bracket  Li-Z dual $\overline{\Psi}_{\nu}$. Then
\[
\sigma f_2\Psi_{\nu}=e^{2 i\pi \nu \beta_2}\Psi_{\nu},
\] 
from \eqref{eq:fgi}. Furthermore,
\begin{eqnarray*}
 Z^{(1)}_{V_{\Z}} \begin{bmatrix}f_1\\ g_1 \end{bmatrix} \left(\Psi_{\nu},w;\overline{\Psi}_{\nu},0; q\right) 
&=& \sum_{\mu\in \Z+\alpha_1} e^{\tpi \mu\beta_1} Z^{(1)}_{\mu}\left(\Psi_{\nu},w;\overline{\Psi}_{\nu},0; q\right). 
\end{eqnarray*}
Thus altogether we find
\begin{eqnarray*}
Z^{(2)}_{V_{\Z}}\begin{bmatrix}f\\g \end{bmatrix}(\tau, w, \rho)
&=&
\sum_{\mu\in \Z+\alpha_1}
\sum_{\nu\in \Z+\kappa}
e^{\tpi \mu\beta_1}
 e^{2 i\pi \nu \beta_2}
\sum_{\Psi_{\nu}\in M\otimes e^{\nu}}
Z^{(1)}_{\mu}\left(\Psi_{\nu},w;\overline{\Psi}_{\nu},0; q\right)
\\
&=&
\sum_{\mu\in \Z+\alpha_1}
\sum_{\nu\in \Z+\kappa}
e^{i\pi(\mu^2\Omega_{11}+2\mu\nu\Omega_{12}+\nu^2\Omega_{22})}
e^{\tpi (\mu\beta_1+\nu \beta_2)}
Z_{M}^{(2)}(\tau,w,\rho )\\
&=&
\vartheta
 \begin{bmatrix}
\alpha \\ \beta  
\end{bmatrix}
(  \Omega)\;
Z^{(2)}_{M}(\tau, w, \rho).\quad \square
\end{eqnarray*}

Comparing with the fermionic expression \eqref{Z2_ferm} and noting that 
\begin{equation}
(e^{i\pi B }\rho)^{\half \kappa^2}\,
Z^{(1)}_{V_{\Z}}
 \begin{bmatrix}f_1\\ g_1 \end{bmatrix}\left( 
 e^{\kappa}, w; \; 
e^{-\kappa}, 0; \; q\right)=\left(\frac{e^{i\pi B }\rho}{ K (w, \tau)^2 }\right)^{\half \kappa^2}
\frac{1}{\eta(\tau)}
\vartheta 
\begin{bmatrix}
\alpha_1 \\ \beta_1  
\end{bmatrix}
 \left(\kappa w ,  \tau \right)
,
\label{factor}
\end{equation}
we obtain a genus two analogue
of the classical Jacobi triple product identity for the elliptic theta function (separate to a similar identity shown in \cite{TZ2}) 
as follows:
\begin{theorem}
\label{theorem:Triple}
The genus two Riemann theta function on ${\mathcal D}^{\rho}$ is given by
\begin{eqnarray*}
\vartheta
 \begin{bmatrix}
\alpha \\ \beta  
\end{bmatrix}
(  \Omega) 
=
e^{2 i\pi \beta_2\kappa} 
\left(\frac{e^{i\pi B }\rho}{ K (w, \tau)^2 }\right)^{\half \kappa^2}
\vartheta 
\begin{bmatrix}
\alpha_1 \\ \beta_1  
\end{bmatrix}
 \left(\kappa w ,  \tau \right)
\det \left(I - T \right) \det \left(I - R \right)^{\half}.
\label{eq:}
\end{eqnarray*}
\end{theorem}
Since $R,T\rightarrow 0$ as $\rho\rightarrow 0$ then 
$\vartheta
 \left[\begin{smallmatrix}
\alpha \\ \beta  
\end{smallmatrix}\right]
( \Omega) 
\approx
e^{2 i\pi \beta_2\kappa }\left(\frac{e^{i\pi B }\rho}{ K (w, \tau)^2 }\right)^{\half \kappa^2}
\vartheta 
\left[
\begin{smallmatrix}
\alpha_1 \\ \beta_1  
\end{smallmatrix}
\right]
 \left(\kappa w ,  \tau \right)
$ to leading order in $\rho$  
as discussed in greater detail in \cite{TZ1}.
As suggested in \cite{R}, we also remark that  Fay's Trisecant Identity at genus two follows by comparing the generating function ${\mathcal G}^{(2)}_n\left[\begin{smallmatrix}f\\ g \end{smallmatrix}\right]
(x_1, y_1,  \ldots, x_n, y_{n})$ of \eqref{gen_function_2} to its 
form in terms of a bosonic basis. We do not give the details here.

\subsection{Modular Invariance}
\label{holomor} 
We lastly consider the modular properties of the genus two  
partition \eqref{Z2_ferm} and $n$--point generating \eqref{Z2_def_epsorbi_two_point_1_alpha} 
functions for the rank two fermion VOSA.    
In \cite{TZ1} we proved the modular invariance of the genus two 
Szeg\H{o} kernel $S^{(2)}
\left[\begin{smallmatrix}\theta^{(2)}\\ \phi^{(2)}\end{smallmatrix}\right]
(x, y \vert \tau, w,\rho)$ of \eqref{s_genus_two}   
under the action of a particular subgroup $L\subset \Sp(4,{\mathbb Z})$ 
and verified that \eqref{Szmod} holds.
We also showed that $\det(I-T)$
is modular invariant. 

Let us recall the relevant definitions used to describe $L$ \cite{MT2}. 
Consider the Heisenberg group $\widehat{H} \subset \Sp(4,{\mathbb Z})$ with elements
\begin{equation}
\mu (a,b,c)=\left( 
\begin{array}{cccc}
1 & 0 & 0 & b \\ 
a & 1 & b & c \\ 
0 & 0 & 1 & -a \\ 
0 & 0 & 0 & 1%
\end{array}%
\right).  \label{mudef}
\end{equation}%
$\widehat{H}$ is generated by $A=\mu (1,0,0)$, $B=\mu (0,1,0)$ and $C=\mu (0,0,1)$ with relations 
$[A,B]C^{-2}=[A,C]=[B,C]=1$. 
We also define $\Gamma _{1}\subset \Sp(4,\Z)$ where $\Gamma_{1}\cong  \SL(2,\Z)$ with elements
\begin{equation}
\gamma_1 =\left( 
\begin{array}{cccc}
a_{1} & 0 & b_{1} & 0 \\ 
0 & 1 & 0 & 0 \\ 
c_{1} & 0 & d_{1} & 0 \\ 
0 & 0 & 0 & 1%
\end{array}%
\right),\ \ a_{1}d_{1}-b_{1}c_{1}=1.
\label{gamma}
\end{equation}%
Together these groups generate $L=\widehat{H}\rtimes\Gamma _{1}\subset \Sp(4,\Z)$ with center $Z(L)=\langle C\rangle $ where $J=L/Z(L)\cong 
\Z^{2}\rtimes  \SL(2,\Z)$ is the Jacobi group. 
From Lemma~15 of \cite{MT2} we find that $L$ acts on $x=(\tau, w, \rho)\in \mathcal{D}^{\rho }$ of 
\eqref{Drho} as follows: 
\begin{eqnarray}
\mu(a, b, c).x &=& (\tau, w+2 \pi i a \tau + 2 \pi i b, \rho),
\label{muDrho} 
\\
\gamma_1. x &=& \left( \frac{a_{1} \tau+b_{1}}{c_{1}\tau + d_{1}}, \frac{w}{c_{1}
\tau + d_{1}}, \frac{\rho}{(c_{1} \tau + d_{1})^2} \right).  
\label{gamDrho}
\end{eqnarray}
Note that the Jacobi group $J$ acts properly on $\mathcal{D}^{\rho }$.

The genus
 two partition function
$Z^{(2)}_{V_{\Z}}\left [\begin{smallmatrix}f\\g \end{smallmatrix}\right]( \tau, w, \rho)$
is not single-valued on $\mathcal{D}^{\rho }$ due to the 
$\left(e^{i\pi B }\rho\, K (w, \tau)^{-2} \right)^{\half \kappa^2}$ factor in \eqref{Z2_ferm}. 
In Section~6 of \cite{MT2} a detailed description is given of the logarithmic function
\begin{equation}  \label{ldef}
l(x)=\log \left (-\frac{\rho }{K(\tau ,w)^{2}}\right), \quad x = (\tau, w, \rho) \in 
\mathcal{D}^{\rho}.
\end{equation}
In particular, we describe  a covering space $\mathcal{\widehat{D}}^{\rho }$ with projection $p:\mathcal{\widehat{D}}^{\rho }\rightarrow \mathcal{{D}}^{\rho }$ where 
$p(\xhat)=(\tau ,w,\rho )$ on which
 $l$ lifts to a single-valued holomorphic function $\widehat{l}:\mathcal{\widehat{D}}^{\rho }\rightarrow \C$. Furthermore, $L$ acts on $\mathcal{\widehat{D}}^{\rho }$   with $p(\gamma.\xhat)=\gamma .x$ for all $\gamma\in L$ where \cite{MT2}
\begin{eqnarray}
\widehat{l}(\mu (a,b,c).\xhat) &=&\widehat{l}(\xhat)+2\pi ia^{2}\tau +2aw
+2\pi i(ab+c), 
\label{mul} 
\\
\widehat{l}(\gamma _{1}.\xhat) &=&\widehat{l}(\xhat)-\frac{1}{2\pi i}
\frac{c_{1}w^{2}}{c_{1}\tau +d_{1}}
\label{gaml}.
\end{eqnarray}%
Thus, we may lift $Z^{(2)}_{V_{\Z}}\left [\begin{smallmatrix}f\\g \end{smallmatrix}\right]( \tau, w, \rho)$ to a holomorphic function 
$\widehat{Z}^{(2)}_{V_{\Z}}\left [\begin{smallmatrix}f\\g \end{smallmatrix}\right](\xhat)$
on $\mathcal{\widehat{D}}^{\rho }$ where the factor $\left(e^{i\pi B }\rho\, K (w, \tau)^{-2} \right)^{\half \kappa^2}$ is lifted to  
$\exp\left(\half\kappa^2 \ \widehat{l}(\xhat)\right)$ and the remaining  single-valued terms of $Z^{(2)}_{V_{\Z}}\left [\begin{smallmatrix}f\\g \end{smallmatrix}\right]\left( \tau, w, \rho\right)$   lift trivially.

We next define an action of $L$ on the automorphisms of   \eqref{eq:fgi}
\begin{eqnarray}
\mu(a,b,c).
\begin{bmatrix}
f_{1}\\ f_{2} \\ g_{1} \\ g_{2}
\end{bmatrix}
=
\begin{bmatrix}
\sigma^b f_{1} g_{2}^b \\ 
\sigma^{a+b+ab} f_{1}^a f_{2} g_{1}^b g_{2}^c  \\ 
\sigma^{-a} g_{1}g_{2}^{-a} \\ g_{2}
\end{bmatrix}
,\quad
\gamma_1.
\begin{bmatrix}
f_{1}\\ f_{2} \\ g_{1} \\ g_{2}
\end{bmatrix}
=
\begin{bmatrix}
f_{1}^{a_1}g_{1}^{b_1}\\ f_{2} \\ f_{1}^{c_1}g_{1}^{d_1} \\ g_{2}
\end{bmatrix}.  
\label{eq:mugam1fg}
\end{eqnarray}
This action is compatible with \eqref{albetatilde} and it is straightforward to confirm that the defining relations for the generators $A,B,C$ are satisfied. We may also determine the $L$ action on the multipliers $\theta_{i},\phi_{i}$ of \eqref{eq:thetaphi} merely by replacing $f_{i}, g_{i}$ in \eqref{eq:mugam1fg}  by $\theta_{i}^{-1},\phi_{i}^{-1}$, respectively, and $\sigma$ by $-1$.\footnote{We note that the action quoted for the $C$ generator in eqn~(112) of \cite{TZ1} contains an erroneous minus sign}
We thus define a canonical  action  of $\gamma \in L$ on the genus two partition function as follows
\begin{eqnarray}
\left. \widehat{Z}^{(2)}_{V_{\Z}}\begin{bmatrix}f\\g \end{bmatrix}\right\vert_{\gamma} (\xhat)
= 
\widehat{Z}^{(2)}_{V_{\Z}}\left( \gamma.\begin{bmatrix}f\\g \end{bmatrix} \right) \left(\gamma.\xhat \right).   
\label{ZV_gamma_1}
\end{eqnarray}

In \cite{TZ1} we showed that $\det(I-T)$ is invariant under the modular group $L$ (with $L$ acting on the multipliers $\theta_{i},\phi_{i}$ as in \eqref{eq:mugam1fg}) so we need only prove modularity for the following  partition function factor 
\[
F
\begin{bmatrix}f\\g \end{bmatrix}(\xhat)
=e^{2 i\pi \beta_2\kappa }\exp\left(\half\kappa^2 \ \widehat{l}(\xhat)\right)
\frac{1}{\eta(\tau)}
\vartheta 
\begin{bmatrix}
\alpha_1 \\ \beta_1  
\end{bmatrix}
 \left(\kappa w ,  \tau \right).
\]
We first analyse the action of  the $\Gamma_1=\SL(2,\Z)$ generators $S=\left(\begin{smallmatrix}0&1\\ -1&0 \end{smallmatrix} \right)$ and $T=\left(\begin{smallmatrix}1&1\\0&1 \end{smallmatrix} \right)$.  We note that
\begin{eqnarray*}
\vartheta 
\begin{bmatrix}
-\beta_{1} \\ 
\alpha_{1}%
\end{bmatrix}
 (\kappa(-\frac{w}{\tau }),-\frac{1}{\tau }) &=&
(-i\tau )^{\half}
e^{-2\pi i\alpha_{1} \beta_{1}}e^{-i\kappa^{2}w^2/4\pi \tau }
\vartheta 
\begin{bmatrix}
\alpha_{1} \\ 
\beta_{1}
\end{bmatrix}
 (\kappa w,\tau )
,  \\
\vartheta 
\begin{bmatrix}
\alpha_{1} \\ 
\beta_{1}-\alpha_{1}-\half%
\end{bmatrix}%
(\kappa w,\tau +1) &=&
e^{-\pi i(\alpha_{1}^2+\alpha_{1})}
\vartheta 
\begin{bmatrix}
\alpha_{1} \\ 
\beta_{1}%
\end{bmatrix}
 (\kappa w,\tau ). 
\end{eqnarray*} 
Employing  \eqref{gaml} it follows that 
\begin{eqnarray*}
F\left.\begin{bmatrix}f\\g \end{bmatrix}\right\vert_{S}(\xhat)
=
\chi\begin{bmatrix}f\\g \end{bmatrix}(S)
F\begin{bmatrix}f\\g \end{bmatrix}(\xhat),  
 \quad
F\left.\begin{bmatrix}f\\g \end{bmatrix}\right\vert_{ T}(\xhat)
=
\chi\begin{bmatrix}f\\g \end{bmatrix}(T)
F\begin{bmatrix}f\\g \end{bmatrix}(\xhat),
\end{eqnarray*}
for $S,T$ multipliers 
\[
\chi\begin{bmatrix}f\\g \end{bmatrix}(S)=
e^{-2\pi i\alpha_1 \beta_1},
\quad
\chi\begin{bmatrix}f\\g \end{bmatrix}(T)=
e^{-\pi i(\alpha_1^2+\alpha_1+\frac{1}{12})}.
\] 
Thus it follows that 
$\widehat{Z}^{(2)}_{V_{\Z}}\left [\begin{smallmatrix}f\\g \end{smallmatrix}\right](\xhat)$ is invariant under all $\gamma_1\in\Gamma_1$ up to a multiplier $\chi\left [\begin{smallmatrix}f\\g \end{smallmatrix}\right](\gamma_1)$ generated by $S,T$.

We next consider the action of the Heisenberg group $\widehat{H} $. From the definition of the theta series and using \eqref{mul} we obtain
\begin{eqnarray*}
&&\vartheta 
\begin{bmatrix}
\alpha_{1}-a \kappa \\ 
\beta_{1}-b \kappa
\end{bmatrix}
 (\kappa( w+2 \pi i a \tau + 2 \pi i b),{\tau }) =
e^{-(\pi i\tau{a}^{2}+aw){\kappa}^{2}}
e^{-2 \pi ia \kappa\beta_{{1}}}
\vartheta 
\begin{bmatrix}
\alpha_{1} \\ 
\beta_{1}
\end{bmatrix}
 (\kappa w,\tau ),
\\
&&\exp\left(\half\kappa^2 \ \widehat{l}(\mu(a,b,c).\xhat)\right)=
e^{(\pi i\tau{a}^{2}+aw+\pi iab+\pi i c){\kappa}^{2}}
\exp\left(\half\kappa^2 \ \widehat{l}(\xhat)\right). 
\end{eqnarray*}  
Furthermore, from \eqref{eq:mugam1fg} we note that 
\[
\mu(a,b,c).\beta_{2}=\beta_{2}+a\beta_{1}-b\alpha_{1}-c\kappa+\half(ab-c).
\]
These results imply that
\[
F\left.\begin{bmatrix}f\\g \end{bmatrix}\right\vert_{ \mu(a,b,c)}(\xhat)
=
\chi\begin{bmatrix}f\\g \end{bmatrix}({\mu(a,b,c)}) F\begin{bmatrix}f\\g \end{bmatrix}(\xhat),
\]  
for multiplier 
\[
\chi\begin{bmatrix}f\\g \end{bmatrix}({\mu(a,b,c)})=
e^{-2\pi i b\alpha_{1}\kappa}
e^{\pi i(ab-c)\kappa(\kappa+1)}.
\] 
Thus, altogether we have shown that 
\begin{theorem}
\label{Theorem_modul_inv}
$\widehat{Z}^{(2)}_{V_{\Z}}\left [\begin{smallmatrix}f\\g \end{smallmatrix}\right](\xhat)$
is holomorphic on $\mathcal{\widehat{D}}^{\rho }$ and  is modular invariant with respect to $\gamma\in L$ with a multiplier system: 
\begin{equation*}
\left. \widehat{Z}^{(2)}_{V_{\Z}}\begin{bmatrix}f\\g \end{bmatrix}\right\vert_{\gamma} (\xhat)
=
\chi\begin{bmatrix}f\\g \end{bmatrix}(\gamma)
\widehat{Z}^{(2)}_{V_{\Z}}\begin{bmatrix}f\\g \end{bmatrix}(\xhat).  
\end{equation*}
\end{theorem}

These modular properties and multiplier system can also be confirmed by comparison to the modular properties of the genus two Riemann theta function expression of Theorem~\ref{Theorem_Z2_fermboson} using results of \cite{MT2}. In particular,  the multiplier $\chi\left [\begin{smallmatrix}f\\g \end{smallmatrix}\right](\mu(a,b,c))$ arises from the genus two theta function and can be easily independently verified.

Finally, in \cite{TZ2} we showed  that the Szeg\H{o} kernel is  modular invariant under $L$ so that we conclude that the genus generating function
 ${\mathcal G}^{(2)}_n \left [\begin{smallmatrix}f\\g \end{smallmatrix}\right]$ lifts to a
 modular invariant form on $\mathcal{\widehat{D}}^{\rho }$ with the same multiplier.

\section{Appendix}

\subsection{Proof of Proposition \ref{prop_assoc}}
\label{PropGen}
To prove \eqref{genloc} we consider the coefficient of $z_0^{-N-1+\alpha \beta}$  of the Jacobi identity \eqref{Jac}. Using \eqref{delta} and \eqref{xylambda} we have
\begin{eqnarray*}
z_0^{-1} 
\left( \frac{z_1 - z_2}{z_0}\right)^{-\alpha \beta}  
\delta\left( \frac{z_1 - z_2}{z_0}\right) &=&
\sum_{r \in \Z}(z_1-z_2)^{r-\alpha \beta} z_0^{-r-1+\alpha \beta}.
\end{eqnarray*}
Thus the coefficient of $z_0^{-N-1+\alpha \beta}$ arising from the first term on the left hand side of \eqref{Jac} 
is 
\[
\left(z_1-z_2\right)^{N-\alpha \beta}
\; {\cal Y}  \left(u\otimes e^{\alpha}, z_1 \right) 
{\cal Y}  \left(v \otimes e^{\beta}, z_2\right).
\]
A corresponding term arises from the second term on the left hand side of \eqref{Jac}. On the other hand, by lower truncation  \eqref{lowertrun}, only finitely many negative integer powers of $z_0$ occur in
\[ 
z_0^{-\alpha\beta}{\cal Y} \left(u \otimes e^{\alpha}, z_0\right)
 (v\otimes e^{\beta}),
\]
so that \eqref{genloc} holds for $N\gg 0$. 
 
Apply the operators of \eqref{genloc} to the vacuum vector together with translation and creativity to find 
\begin{eqnarray*}
&&
\left(z_1-z_2\right)^{N}\left(z_1-z_2\right)^{-\alpha \beta} e^{z_2L(-1)}
{\cal Y}  \left(u\otimes e^{\alpha}, z_1 -z_2\right) 
v \otimes e^{\beta}\\
&=&
\left(z_1-z_2\right)^{N}\left(z_2-z_1\right)^{-\alpha \beta}  e^{z_1 L(-1)}
{\cal Y}  \left(v \otimes e^{\beta}, z_2-z_1\right)
u\otimes e^{\alpha}.
\end{eqnarray*}
Hence \eqref{genskew} follows for $z=z_1-z_2$.

To prove \eqref{genassoc} we consider the action on $w\otimes e^{\gamma}$ of the vertex operators appearing in the generalized Jacobi identity \eqref{Jac}.  The first term on the left hand side of \eqref{Jac} gives 
\begin{eqnarray}
{\cal Y}  \left(u\otimes e^{\alpha}, z_1 \right)  
{\cal Y}  \left(v \otimes e^{\beta}, z_2\right) w\otimes e^{\gamma}
=
z_1^{\alpha( \beta+\gamma)}z_2^{\beta \gamma}
A(z_1)B(z_2) w\otimes e^{\alpha+\beta+\gamma},\quad
\label{Jac1}
\end{eqnarray}
for $A(z_1)=Y_{-}(\alpha, z_1) Y(u ,z_1)Y_{+}(\alpha,z_1)$ and $ B(z_2)=Y_{-}(\beta, z_2) Y(v,z_2)Y_{+}(\beta,z_2)$.
We note that $A(z_1)B(z_2) w\otimes e^{\alpha+\beta+\gamma}=\sum_{k\in \Z}A(k)z_1^k B(z_2) w\otimes e^{\alpha+\beta+\gamma}$ for some operators $A(k)$. Furthermore
\begin{eqnarray*}
z_0^{-1} 
\left( \frac{z_1 - z_2}{z_0}\right)^{-\alpha \beta}  
\delta\left( \frac{z_1 - z_2}{z_0}\right)  z_1^{k+\alpha \beta}&=&
\sum_{r \in \Z}\sum_{s\ge 0}\binom{r-\alpha \beta}{s} (-z_2)^sz_0^{-r+\alpha \beta-1} z_1^{r+k-s}\\
&=&
\sum_{N \in \Z}z_1^{-N-1}(z_0+z_2)^{N+k+\alpha\beta}.
\end{eqnarray*}
Thus the coefficient of $z_1^{-N-1+\alpha \gamma}$  of \eqref{Jac1} is  
\begin{eqnarray*}
&&(z_0+z_2)^{N+\alpha \beta}z_2^{\beta \gamma} A(z_0+z_2)B(z_2) \;w\otimes e^{\alpha+\beta+\gamma}\\
&&=\left(z_0+z_2\right)^{N-\alpha \gamma}
\; {\cal Y}  \left(u\otimes e^{\alpha}, z_0+z_2 \right) 
{\cal Y}  \left(v \otimes e^{\beta}, z_2\right) \,  w\otimes e^{ \gamma}.
\end{eqnarray*}
On the other hand, the second term on left hand side of the Jacobi identity gives
\begin{eqnarray*}
&&
- z_0^{-1}
\left( \frac{z_2 - z_1}{z_0}\right)^{ -\alpha\beta }  
 \delta\left( \frac{z_2 - z_1}{-z_0}\right)
 \;{\cal Y} \left(v \otimes e^{\beta}, z_2\right)   {\cal Y} \left(u\otimes e^{\alpha}, z_1 \right)
w\otimes e^{\gamma}
\\
&&
=-z_0^{-1} 
\left( \frac{z_2 - z_1}{z_0}\right)^{-\alpha \beta}  
\delta\left( \frac{z_2 - z_1}{-z_0}\right) z_1^{\alpha\gamma} z_2^{ \beta(\alpha+\gamma)}
\; B(z_2) A(z_1)\;w\otimes e^{\alpha+\beta+\gamma}.
\end{eqnarray*}
The coefficient of $z_1^{-N-1+\alpha \gamma}$ of this expression vanishes for $N\gg0$ by lower truncation.

The right hand side of the Jacobi relation gives
\begin{eqnarray*}
&&
z_2^{-1}    
\delta\left( \frac{z_1 - z_0}{z_2}\right)
{\cal Y} \left( {\cal Y} \left(u \otimes e^{\alpha}, z_0\right)
 (v\otimes e^{\beta}), z_2\right)
\left( \frac{z_1 - z_0}{z_2}\right)^{\alpha a(0)}w\otimes e^{\gamma}
\\
&&=
z_2^{-1}    
\left( \frac{z_1 - z_0}{z_2}\right)^{\alpha \gamma}\delta\left( \frac{z_1 - z_0}{z_2}\right)
{\cal Y} \left( {\cal Y} \left(u \otimes e^{\alpha}, z_0\right)
 (v\otimes e^{\beta}), z_2\right)\;w\otimes e^{\gamma}.
\end{eqnarray*}
But 
\begin{eqnarray*}
z_2^{-1}    
\left( \frac{z_1 - z_0}{z_2}\right)^{\alpha \gamma}\delta\left( \frac{z_1 - z_0}{z_2}\right)
&=&
\sum_{r \in \Z}\sum_{s\ge 0}\binom{r+\alpha \gamma}{s} (-z_0)^sz_1 ^{r-s+\alpha \gamma} z_2^{-r-\alpha \gamma-1}\\
&=&
\sum_{N \in \Z}z_1^{-N-1+\alpha \gamma}(z_2+z_0)^{N-\alpha\gamma}.
\end{eqnarray*}
Taking the coefficient of $z_1^{-N-1+\alpha \gamma}$ of this expression leads to the result. 
Finally, \eqref{genComm} follows in the usual way for a VOA an substituting $\alpha=0$ into \eqref{Jac} $\square$
\medskip

\subsection{Proof of Proposition~\ref{prop_ZalphaK}}
\label{ProofZalphaK}
We first show that the general  intertwining  $n$--point function can be written in terms of a 1-point function:
\begin{eqnarray}
Z^{(1)}_{\alpha}  
\left( e^{\beta_1}, z_1;  \ldots ;e^{\beta_n}, z_n;q \right)
&=&Z^{(1)}_{\alpha}\left( 
\Ycal \left[ e^{\beta_1}, z_{1n}\right]
\ldots\Ycal \left[e^{\beta_{n-1}}, z_{n-1n}\right]
e^{\beta_n};q
\right)\qquad\qquad
\label{eq:Znto1a}
\\
&=&Z^{(1)}_{\alpha}\left(\Ycal[e^{\beta_1}, z_1]  \ldots \Ycal[e^{\beta_n},z_n]\vac;q\right),
\label{eq:Znto1b}
\end{eqnarray}
where $\Ycal[e^{\beta_1}, z_1]  \ldots \Ycal[e^{\beta_n},z_n]\vac\in M[[z_i,z_i^{-1}]]$ since  $\sum_{i=1}^{n}\beta_i= 0$.
Using \eqref{pcom1}, \eqref{eq:Lgn} and \eqref{eq:qhatconj} with $\Delta(\alpha,q_i)e^{\beta_i}=q_i^{\alpha a(0)}e^{\beta_i}$ we find
\begin{eqnarray*}
&&Z^{(1)}_{\alpha}  
\left( e^{\beta_1}, z_1;  \ldots ;e^{\beta_n}, z_n;q \right)\\
&&=
\Tr_{M}\left(
e^{-\alpha \qhat } 
\Ycal \left(q_{1}^{L(0)} e^{\beta_1}, q_{1}\right) 
\ldots 
\Ycal \left(q_{n}^{L(0)}e^{\beta_n}, q_{n}\right) 
q^{L(0)-1/24} 
e^{\alpha \qhat }
\right)\\
&&=
\Tr_{M}\left( 
\Ycal \left(q_{1}^{L(0)+\alpha a(0)}e^{\beta_1}, q_{1}\right) 
\ldots 
\Ycal \left(q_{n}^{L(0)+\alpha a(0)} e^{\beta_n}, q_{n}\right) 
q^{L_{g}(0)-1/24} 
\right).
\end{eqnarray*}
Define the ``shifted'' conformal vector  
$
\omega_{\alpha}=\omega-\alpha a(-2)\vac
$
with modes \cite{DM}
\[
L_{\alpha}(n)=L(n)+\alpha(n+1)a(n),
\]
which satisfy the Virasoro algebra with central charge $c_\alpha = 1-12 \alpha^2$ so that
 $(\Mcal,\Ycal,\vac,\omega_{\alpha})$ is a Generalized VOA of central charge $c_\alpha$. 
Noting  that
\[
L_g(0)-\frac{1}{24}=L_{\alpha}(0)-\frac{c_\alpha}{24}, 
\]
we thus obtain (cf.  Proposition~9 of \cite{MTZ})
\begin{eqnarray*}
Z^{(1)}_{\alpha}  
\left( e^{\beta_1}, z_1;  \ldots ;e^{\beta_n}, z_n ;q\right)
=\Tr_{M}\left( 
\Ycal \left(q_{1}^{L_{\alpha}(0)}e^{\beta_1}, q_{1}\right) 
\ldots 
\Ycal \left(q_{n}^{L_{\alpha}(0)} e^{\beta_n}, q_{n}\right) 
q^{L_{\alpha}(0)-c_{\alpha}/24} 
\right).
\end{eqnarray*}
Since we are tracing over $M$, we may consecutively apply the standard form of associativity   to obtain:
\begin{eqnarray}
&&\Tr_{M}
\left( 
Y \left(
\Ycal \left(q_{1}^{L_{\alpha}(0)} e^{\beta_1}, q_1-q_n\right) 
\Ycal \left(q_{2}^{L_{\alpha}(0)}e^{\beta_2}, q_2-q_n\right) 
\ldots
\right.\right.
\notag
\\
&&
\left.\left.
\ldots\Ycal \left(q_{n-1}^{L_{\alpha}(0)}e^{\beta_{n-1}}, q_n-q_n\right) 
q_{n}^{L_{\alpha}(0)} e^{\beta_n}, q_{n}\right) 
q^{L_{\alpha}(0)-c_{\alpha}/24} 
\right)
\notag
\\
&&=\Tr_{M}\left( 
Y\left (q_{n}^{L_{\alpha}(0)}
\Ycal \left[ e^{\beta_1}, z_{1n}\right]_{\alpha} 
\ldots\Ycal \left[e^{\beta_{n-1}}, z_{n-1n}\right]_{\alpha} 
e^{\beta_n}, q_{n}\right) 
q^{L_{\alpha}(0)-c_{\alpha}/24} 
\right),\qquad\quad
\label{Zn1alpha}
\end{eqnarray}
on using 
\[
q_{n}^{-L_{\alpha}(0)}\Ycal(e^{\beta},q_{i})q_{n}^{L_{\alpha}(0)}
=\Ycal(q_{n}^{-L_{\alpha}(0)}e^{\beta},q_{n}^{-1}q_{i}),
\]
 and where
\[
\Ycal \left[e^{\beta}, z \right]_{\alpha}=\Ycal \left(q_z^{L_{\alpha}(0)}e^{\beta}, q_z-1\right),
\]
denotes the square bracket vertex operator for the shifted Virasoro system. 
But  \eqref{eq:qhatconj} together with the identities
\begin{eqnarray*}
\Delta(-\alpha,q_n)
q_{n}^{L_{\alpha}(0)}
&=&
q_{n}^{L(0)}\Delta(-\alpha,1), \\
\Delta(-\alpha,1)\Ycal \left[e^{\beta}, z \right]_{\alpha}\Delta(\alpha,1)
&=&
\Ycal \left[e^{\beta}, z \right],
\end{eqnarray*} 
imply that \eqref{Zn1alpha} becomes the 1-point function \eqref{eq:Znto1a}.
\eqref{eq:Znto1b} follows since
\begin{eqnarray*}
o(\Ycal \left[ e^{\beta_1}, z_{1}\right]
\ldots
\Ycal \left[e^{\beta_{n}}, z_{n}\right]\vac
)
&=&
o(\Ycal \left[ e^{\beta_1}, z_{1}\right]
\ldots\Ycal \left[e^{\beta_{n-1}}, z_{n-1}\right] e^{z_nL[-1]}e^{\beta_{n}}
)\\
&=&
o(e^{z_nL[-1]}\Ycal \left[ e^{\beta_1}, z_{1n}\right]
\ldots\Ycal \left[e^{\beta_{n-1}}, z_{n-1n}\right] e^{\beta_{n}}
)\\
&=&
o(\Ycal \left[ e^{\beta_1}, z_{1n}\right]
\ldots\Ycal \left[e^{\beta_{n-1}}, z_{n-1n}\right] e^{\beta_{n}}
),
\end{eqnarray*}
and since $o(L[-1]u)=0$ for all $u\in M$ \cite{Z1}.  

As in  the proof of  Proposition~5 of \cite{MT1}, by applying the identity  \eqref{Ypmcom} we find that \eqref{eq:Znto1b} is given by
\[
\prod_{1\leq r<s\leq n}z_{rs}^{\beta _{r}\beta _{s}}
Z^{(1)}_{\alpha}
\left(
\exp \left(
\sum_{m>0}\frac{a[-m]}{m}\sum_{i=1}^{n}\beta _{i}z_{i}^{m}
\right)\vac;q
\right).  
\]
Proposition~3 of \cite{MT1} gives 
\begin{eqnarray*}
&&Z^{(1)}_{\alpha}
\left(
\exp \left(
\sum_{m>0}\frac{a[-m]}{m}\sum_{i=1}^{n}\beta _{i}z_{i}^{m}
\right)\vac;q
\right) 
\\
&=&
q^{\half \alpha^2}
\exp\left(\alpha \sum_{i=1}^{n}\beta_{i}z_{i}\right)
Z_{0}\left(\exp (\sum_{m>0}\frac{a[-m]}{m}\sum_{i=1}^{n}\beta _{i}
z_{i}^{m})\vac;q \right).
\end{eqnarray*}
Finally, Proposition~4 of \cite{MT1} gives
\begin{equation*}
Z_{0}(\exp (\sum_{m>0}\frac{a[-m]}{m}\sum_{i=1}^{n}\beta _{i}z_{i}^{m})\vac;q)
=\frac{1}{\eta (q )}
\prod_{1\leq r<s\leq n}\left(\frac{K(z_{rs},\tau )}{z_{rs}}\right)^{\beta _{r}\beta _{s}}. 
\end{equation*}
The proposition follows. $\square$


\end{document}